\documentclass{article}%
\usepackage{amsfonts}
\usepackage{graphicx}
\usepackage{amsmath}
\usepackage{amssymb}%
\providecommand{\U}[1]{\protect\rule{.1in}{.1in}}
\newcommand{\Monm}{\mathcal{F}^*}
\newcommand{\Mon}{{\rm Mon}^*}
\newcommand{\tMonm}{\mathcal{F}}
\newcommand{\tMon}{{\rm Mon}}

\newcommand\R{{\mathbf {R}}}
\newcommand\Z{{\mathbf {Z}}}

\newcommand\bkE{{\mathbf {E}}}
\newcommand\E{{\mathbf {E}}}

\newtheorem {Lemma}{Lemma}[section]

\newtheorem {Theorem}{Theorem}[section]
\newtheorem {Proposition}{Proposition}[section]
\newtheorem {Corollary}{Corollary}[section]

\newtheorem{Definition}{Definition}[section]

\newtheorem{Remark}{Remark}[section]

\newcommand\I{{ 1\hspace{-1,2mm}{\mathrm I}}}

\newcommand\F{{\mathcal {F}}}

\newcommand\beq{\begin{equation}}
\newcommand\eeq{\end{equation}}

\begin{document}
\begin{center} {\bf \Large Rates of convergence  in
the strong invariance principle under projective criteria}\vskip15pt

J\'er\^ome Dedecker $^{a}$, Paul Doukhan, $^{b}$ {\it
and\/} Florence Merlev\`{e}de   $^{c}$
\end{center}
\vskip10pt
$^a$ Universit\'e Paris 5. Email: jerome.dedecker@parisdescartes.fr\\ \\
$^b$ Universit\'e Cergy Pontoise. E-mail: Doukhan@u-cergy.fr\\ \\
$^c$ Universit\'e Paris Est-Marne La Vall\'ee.
E-mail: florence.merlevede@univ-mlv.fr \vskip10pt
{\it Key words}: almost sure invariance principle, strong approximations, weak dependence,   Markov chains.\vskip5pt

\begin{center}
{\bf Abstract}\vskip10pt
\end{center}

We give rates of convergence in the strong invariance principle
for stationary  sequences satisfying some projective
criteria. The conditions are expressed in terms  of conditional expectations of partial sums of the initial sequence. Our results apply to a large variety of examples. We present some  applications to a reversible Markov chain, to symmetric random walks on the circle,
and to  functions of dependent sequences.

\bigskip

\section{Introduction and notations} \label{intronota}
\setcounter{equation}{0}
The almost sure invariance principle is a powerful tool in both  probability and statistics. It says that the partial sums of  random variables can be approximated by those of independent Gaussian random variables, and that the approximation error between the trajectories of the two processes is negligible in a certain sense. In this paper, we are interested in studying rates in the almost sure invariance principle for dependent sequences.

When $(X_i)_{i \geq 1}$ is a sequence of independent and identically distributed (iid) centered real-valued random variables with a finite second moment, it is known from Strassen (1964) that a sequence $(Z_i)_{i \geq 1}$ of iid centered  Gaussian variables with variance $\sigma^2 = \E (X_0^2)$ may be constructed is such a way that
\beq \label{strassen}
\sup_{1 \leq k \leq n} \Big| \sum_{i=1}^k (X_i - Z_i ) \Big| = o(b_n) \text{ almost surely, as $n\rightarrow \infty$},
\eeq
where $b_n = (n\log \log n)^{1/2}$. To get smaller $(b_n)$ additional information on the moments of $X_1$ is necessary. In the iid setting, Koml\'os, Major and Tusn\'ady (1976) and  Major (1976)  obtained (\ref{strassen}) with $b_n = n^{1/p}$ as soon as ${\mathbf E}(|X_1|^p) < \infty$  for $p>2$.

There has been a great amount of works to extend these results to dependent sequences, under various  conditions: see for instance Heyde (1975), Philipp and Stout (1975), Berkes and Philipp (1979), Dabrowski (1982), Bradley (1983), Utev (1984), Eberlein (1986),  Shao and Lu (1987), Sakhanenko (1988),
Shao (1993), Rio (1995), and more recently, Wu (2007), Zhao and Woodroofe (2008),
Liu and Lin (2009),  Gou\"ezel (2010), Merlev\`ede and Rio (2012).

Having explicite rates in the strong invariance principle \eqref{strassen}
may be useful to  derive results in asymptotic statistics. We refer to the monograph by Cs\"org\H{o} and Horv\'ath (1997) which illustrates the importance of strong approximation principles for change-point and trend analysis, and also to the paper by Horv\'ath and Steinebach (2000), showing  that limit results for so-called CUSUM and MOSUM-type test procedures, which are used to detect mean and variance changes, can be proved with the help of strong invariance principles. For instance, Aue, Berkes and Horv\'ath (2006) use a strong approximation principle with an explicit rate for the  sums of  squares of augmented GARCH sequences to study the limiting behavior of statistical tests, which are used to decide whether the volatility of the underlying variables is stable over time or if it changes in the observation
period. Let us also  mention the recent paper by Wu and Zhao (2007) who consider statistical inference of trends in mean non-stationary models.
Starting from a strong approximation principle with an explicit rate for the partial sums of stationary processes, they propose a statistical test
concerning  the existence of structural breaks in trends, and
they construct simultaneous confidence bands with
asymptotically correct nominal coverage probabilities. In their paper, they point out that  an explicit rate in the strong approximation principle is crucial to control certain errors terms (see their Remark 2).

In this paper, we  obtain rates of convergence of order $b_n =n^{1/p}L(n)$ in (\ref{strassen}) ($L(n)$ is a slowly varying function) when $p \in ]2,4]$,  for stationary sequence satisfying some
projective conditions. To describe our results more precisely, we need  to introduce some notations.

Let $(\Omega,{\cal A}, {\mathbf P} )$ be a
probability space, and $T:\Omega \mapsto \Omega$ be
a bijective bimeasurable transformation preserving the probability ${\mathbf P} $.
For a $\sigma$-algebra ${\cal F}_0 $ satisfying ${\cal F}_0
\subseteq T^{-1 }({\cal F}_0)$, we define the nondecreasing
filtration $({\cal F}_i)_{i \in {\mathbf Z}}$ by $ {\cal F}_i=T^{-i
}({\cal F}_0)$.
Let $X_0$ be a square integrable,  zero mean and ${\mathcal F}_0$-measurable random variable, and
define the stationary sequence  $(X_i)_{i \in \mathbf Z}$ by
 $X_i = X_0 \circ T^i$. Define then the partial sum $S_n= X_1+X_2+ \cdots  + X_n$. Finally, let $H_i$ be the space of ${\mathcal
F}_i$-measurable and square integrable random variables, and
denote by $H_i \ominus H_{i-1}$ the orthogonal of $H_{i-1}$ in
$H_i$. Let $P_i$ be the projection operator from ${\mathbf L}^2$
to $H_i \ominus H_{i-1}$, that is
$$
     P_i(f)={\mathbf E}(f|{\mathcal F}_i)-{\mathbf E}(f|{\mathcal F}_{i-1})
     \quad \text{for any $f$ in ${\mathbf L}^2$.}
$$
We shall denote sometimes by ${\mathbf E}_i$ the conditional expectation with respect to ${\mathcal F}_i$. The following notations will be also frequently used: For any two
positive sequences $a_{n}\ll b_{n}$ means that for a certain numerical
constant $C$ not depending on $n$, we have $a_{n}\leq Cb_{n}$ for all $n$;
$[x]$ denotes the largest integer smaller or equal to $x$.

Our starting point is the same as in Shao and Lu (1987) and Wu (2007):
to obtain a rate of order $b_n =n^{1/p}L(n)$ in (\ref{strassen}), we shall
always assume  that
\begin{equation}\label{Gor}
  \text{the series} \quad d_0=\sum_{i \geq 0} P_0(X_i) \quad
  \text{converges in ${\mathbf L}^p$.}
\end{equation}
 We then define the approximating martingale
$M_n$ as in Gordin (1969) and Heyde (1974):
\begin{equation}\label{Gor2}
M_n= \sum_{i=1}^n d_0 \circ T^i \, .
\end{equation}
Now to prove (\ref{strassen}), it remains to find appropriate conditions under which
(\ref{strassen}) is true for $M_n$ instead of $S_n$, and $|M_n-S_n|=o(b_n)$ almost surely. To prove that
$M_n$ satisfies (\ref{strassen}), we shall apply  the Skorohod embedding theorem (see Proposition \ref{propmart} in the appendix).

 Let us now describe the main differences between our paper and that by Shao and Lu (1987)
 or  Wu (2007).

 It is quite easy to see that one of the assumptions of Shao and Lu is that
 the sequence
$ {\mathbf E}(S_n|{\mathcal F}_0)$
 converges in ${\mathbf L}^p$, which is equivalent to a coboundary decomposition:
 $X_0=d_0+ Z- Z \circ T$, for some random variable $Z$ in ${\mathbf L}^p$. Clearly this coboundary decomposition implies (\ref{Gor}). However,
 for $p=2$,
 this condition is known to be  too restrictive for the almost sure invariance
 principle: see the recent paper by Zhao and Woodroofe (2008). We shall not require this
 coboundary decomposition in  our  results.

 In his 2007's paper, Wu does not assume the existence of a coboundary decomposition, but 
 a polynomial decay of 
 $$
 \sum_{i \geq n} \|P_0(X_i)\|_p  .
 $$
He also assumes 
 that the quantity $\|\bkE  (d_n^2 | {\mathcal F}_0  ) -
\bkE  (d_0^2)  \Vert_{p/2}$ converges to zero fast enough as $n$ tend to infinity. He then gives a large class
of functions of iid sequences to which his results apply. In our first result (Theorem \ref{interwunous}), we give slightly weaker conditions than those required in Theorem 4 in Wu (2007), and we provide some examples to which our result applies whereas Wu's conditions are not satisfied.

However the condition on $d_n$ can be very difficult to check if the
sequence $X_i$ has not an explicit expression as a function of an iid sequence.  In  Theorem \ref{ThNAS} and its corollaries, we give conditions expressed  in terms of conditional expectations of the random variables $X_i$ and $X_iX_j$ with respect to the past $\sigma$-algebra ${\mathcal F}_0$, to obtain rates in the almost sure invariance principle. The proofs of these results are postponed to the section \ref{sectproofs}.

As we shall see, our results apply to a large variety of examples, including mixing processes
of different kinds. However, with this direct approximating martingale method,
there seems to be no hope to get the rate $n^{1/p}$ instead of $n^{1/p}L(n)$ with only a moment of order $p$, whereas this can be done in some situations $via$ other approaches (for $\phi$-mixing sequences in the sense of Ibragimov (1962) 
and $2<p<5$,
it can be deduced from a paper by  Utev (1984)).

In Section \ref{sectionappli}, we have chosen to restrict our attention to four different classes of examples.

We first apply Theorem \ref{interwunous} to a function of an absolutely regular Markov
chain (see Section \ref{sectinvert}). We obtain a rate of order $n^{1/p}$ in \eqref{strassen} under the same conditions implying a Rosenthal type  inequality of order $p$ (see Rio (2009)). Note that we get these rates of convergence in the case where the $\beta$-mixing coefficients of the chain are not summable.

For the three other classes of examples, we apply Theorem \ref{ThNAS}.

In Section \ref{sectcircle} we show that our projective conditions apply to the well
known example of the symmetric random walk on the circle. We obtain rates of convergence in the strong invariance principle for a function  $f$ of  the stationary
Markov chain with transition $Kf (x) =
\frac{1}{2}( f (x+a) + f (x-a) )$, when $a$ is irrational and  badly approximable by rationals (see definitions \eqref{badly} and \eqref{badlyweak}), and the Fourier coefficients
$\hat f$ of $f$ satisfy $\hat f(k) =O(k^{-b})$ for some $b>1$. In particular, we obtain the rate $n^{1/4} L(n)$ in \eqref{strassen} when $f$ is three times differentiable (see Remark \ref{rkrw}). Up to our knowledge, this is the first strong approximation result
for this chain.

In Section \ref{secttau}, we give an application of Theorem \ref{ThNAS} to the case of $\tau$-dependent sequences
in the sense of Dedecker and Prieur (2005). The nice coupling properties of
$\tau$-dependent sequences enables to get results for sums of H\"older functions of
the random variables. We apply our results  to a functional auto-regressive
process whose  auto-regression function is not strictly contracting, so that the $\tau$-dependence coefficients decrease with an arithmetical rate.

In Section \ref{sectionappliwms}, we give an application of Theorem \ref{ThNAS} to the case of $\alpha$-dependent sequences
in the sense of Dedecker and Prieur (2005). This class contains the class of $\alpha$-mixing sequences in the sense of Rosenblatt,
and as a consequence, we improve on the result given by Shao and Lu (1987) in the $\alpha$-mixing case.
We also give an example of a non $\alpha$-mixing sequence to which our result apply, by considering
the Markov chain associated to an intermittent map of the interval.





\section{Main results.}

\setcounter{equation}{0}
\label{SectionPC} \setcounter{equation}{0} In this section, we
give rates of convergence in the strong invariance principle for
stationary sequences satisfying projective criteria. We shall use the notations of Section \ref{intronota} (recall in particular that $d_0 = \sum_{i \geq 0} P_0(X_i)$,
where  the series converges in ${\mathbf L}^p$).

\medskip

We start this section by recalling Theorem 4 in Wu (2007).

\begin{Theorem} (Wu (2007)). \label{thm4wu}
Let $2 <p \leq 4$. Assume that
\begin{equation} \label{condWUASIP}
\sum_{\ell \geq n}\Vert P_0(X_{\ell}) \Vert_p = O (n^{-(1/2-1/p)} ) \ \text{ and } \
\Vert  \bkE  (d_n^2 | {\mathcal F}_0  ) -
\bkE  (d_n^2)  \Vert_{p/2} = O(n^{-(1-2/p) } ) \, .
\end{equation}
Then, enlarging $\Omega$ if
necessary, there exists a sequence $(Z_i)_{i
\geq 1}$ of iid gaussian random variables with zero mean and
variance $\sigma^2={\mathbf E}(d_0^2)$ such that 
$$
\sup_{1\leq k \leq n} \Big|S_k - \sum_{i=1}^k Z_i\Big|  = o(n^{1/p} ( \log n)^{3/2}   ) \text{ almost surely, as $n\rightarrow \infty$}.
$$
\end{Theorem}
Let us mention that in the statement of Theorem 4 in Wu (2007), the bound in the first part of condition \eqref{condWUASIP} appears with the power $-(1-2/p)$ (instead of $-(1/2-1/p)$ but his proof reveals that it is a missprint in the statement. Let us continue now with some comments concerning the method used to prove this result.
The second part of condition \eqref{condWUASIP} comes from an application of the Shorohod-Strassen embedding theorem (see Proposition \ref{propmart} given in Appendix for more details) to the martingale $M_n$. Hence the conclusion of Theorem \ref{thm4wu} holds provided that
\begin{equation} \label{convergenceps}
R_n = S_n - M_n= o(n^{1/p} ( \log n)^{3/2}   )\text{ almost surely} \, ,
\end{equation}
which is true if
\begin{equation} \label{condcuny1}
\sum_{n\ge 2}\frac{\| R_n \|_{p }^p}{n^{2} ( \log n)^{p/2}} <\infty \, .
\end{equation}
(see Proposition 1 in Wu (2007) combined with an application of H\"older's inequality). Next, Wu proved the following upper bound (see his Theorem 1):
\begin{equation} \label{condwurn}
\Vert R_n \Vert_p^2 \ll \sum_{k=1}^n
 \Big (  \sum_{\ell \geq k}\Vert P_0(X_{\ell}) \Vert_p \Big )^2 \, ,
 \end{equation}
which leads to the first part of condition (\ref{condWUASIP}).  
In Proposition \ref{propproofth22} of Section \ref{sectproofs},
we give another condition under which   (\ref{convergenceps}) 
is satisfied. As a consequence we obtain   the following  result:
\begin{Theorem} \label{interwunous} Let $2< p \leq 4$ and  $t >2/p$.
Assume that $X_0$ belongs to ${\mathbf L}^p$ and that (\ref{Gor}) is satisfied. Assume in addition that
\begin{equation} \label{newcondTH22}
\sum_{n\ge 2}\frac{\| \E (S_n |{ \mathcal F}_0) \|_{p }^p}{n^2 (\log n)^{(t-1)p/2}} <\infty \ \text{ and } \
\sum_{n\ge 2} \frac{\| \E (M^2_n |{ \mathcal F}_0) -\E (M^2_n) \|_{p/2 }^{p/2}}{n^2 (\log n)^{(t-1)p/2}}   <\infty  \, .
\end{equation}Then $n^{-1}\E(S_n^2)$ converges to  some nonnegative number $\sigma^2$
 and, enlarging $\Omega$ if
necessary, there exists a sequence $(Z_i)_{i
\geq 1}$ of iid gaussian random variables with zero mean and
variance $\sigma^2$ such that
$$
\sup_{1\leq k \leq n} \Big|S_k - \sum_{i=1}^k Z_i\Big|  = o \big ( n^{1/p} ( \log n)^{(t+1)/2}  \big ) \text{ almost surely, as $n\rightarrow \infty$}.
$$
\end{Theorem}

Since  
$\| \E (S_n |{ \mathcal F}_0) \|_{p } =\| \E (R_n |{ \mathcal F}_0) \|_{p } \leq \|R_n \|_{p }  $, by taking into account \eqref{condwurn}, it follows that  the first part of \eqref{newcondTH22} holds under the first part of \eqref{condWUASIP}. Therefore Theorem \ref{interwunous} contains Theorem \ref{thm4wu}.

Notice also that the first part of \eqref{newcondTH22} can be satisfied whereas the first part of \eqref{condWUASIP} fails to hold. Indeed, let us consider the following linear process $(X_k)_{ k \in \Z}$ defined by $X_k = \sum_{j \geq 0} a_j \varepsilon_{k-j}$ where $(\varepsilon_k)_{ k \in \Z}$ is a strictly stationary sequence of martingale differences in ${\mathbf L}^p$ and $(a_k)_{ k \in \Z}$ is a sequence of reals defined by:
$$
a_0=1+u_0 \ \text{ and } \ a_k = \frac{1}{k^{\alpha +1}} + (-1)^k u_k \text{ for all $k \geq 1$} \, ,
$$
where $\alpha >0$ and $(u_k)_{ k \in \Z}$ is a sequence of reals in $\ell^2$ but not in $\ell^1$. Taking ${\mathcal F}_0 = \sigma(\varepsilon_k , k \leq 0)$, it follows that $P_0(X_i) =a_i \varepsilon_0 $, showing that (\ref{Gor}) is satisfied but the first part of \eqref{condWUASIP} is not. In addition, from Burkholder's inequality, $$ \| \E (S_n |{ \mathcal F}_0) \|_{p }^2 \ll \Vert \varepsilon_0 \Vert^{2}_p \sum_{j \geq 0} \Big (\sum_{k=j+1}^{n+j}a_k \Big)^2 \ll n^{1-2\alpha} + \sum_{k \geq 0} u^2_k \, .$$ Hence the first part of \eqref{newcondTH22} is satisfied as soon as $\alpha \geq 1/2 -1/p$. In this situation, notice that the second part of \eqref{newcondTH22} is satisfied as soon as it is with $\sum_{k=1}^n \varepsilon_k$ instead of $M_n$. This last condition for $\sum_{k=1}^n \varepsilon_k$ can be then verified in different situations. We refer, for instance, to Section 4.3 in Merlev\`ede and Peligrad (2012) where this condition has been verified in case when the sequence of martingale differences, $(\varepsilon_k)_{ k \in \Z}$, has an ARCH($\infty$) structure, or also to the example given in Section 4.1 in Dedecker {\it et al.} (2009) where  $(\varepsilon_k)_{ k \in \Z}$ is, in addition, a certain function of a homogeneous Markov chain as described in Davydov (1973).

\medskip

Theorem \ref{thm4wu} as well as Theorem \ref{interwunous}  gives  explicit approximation rates that are optimal up to multiplicative logarithmic factors. As we just mentioned before, the conditions involved in these results are well adapted to linear processes even generated by martingale differences sequences, and we would like  to refer to Section 3 in Wu (2007) where it is shown that they are also satisfied for a large variety of functions of iid sequences. In Section \ref{sectinvert}, we shall also give an application of Theorem \ref{interwunous} to the case where $(X_n)_{n \geq 0}$ is a function of a stationary Markov chain for which the knowledge of the transition probability allows us to verify  both parts of condition \eqref{newcondTH22}.

However, a  condition expressed in terms of
$\Vert \E(S_n^2 |{\mathcal F}_0) -  \E(S_n^2) \Vert_{p/2}$
rather than the second part of condition \eqref{newcondTH22} would give  a nice counterpart
to  Theorem
\ref{interwunous}.
It would be much easier to check, and  would allow to consider general classes of weakly dependent processes that are not explicit functions of iid sequences. The forthcoming Theorem \ref{ThNAS} and its corollaries are in this direction.

To replace $M_n$ by $S_n$ in the second part of condition \eqref{newcondTH22}, a first step is to  give a  precise decomposition of $R_n=S_n-M_n$.

\begin{Proposition} \label{approxrnq} Let  $p \geq 1$ and assume that \eqref{Gor} holds. Then, for any positive integers $n$ and  $N$,
\begin{enumerate}
\item $
R_n  =  \E (S_n |{\mathcal F}_0) -\bkE (S_{n+N} - S_n |{\mathcal F}_{n}) + \bkE (S_{n+N} - S_n  |{\mathcal F}_{0}) -  \sum_{k=1}^n \sum_{j \geq n+N +1} P_k(X_j)
 \, .
$
\item  $
\Vert R_n \Vert^{p'}_{p} \ll   \Vert \bkE (S_{n}   |{\mathcal F}_{0}) \Vert_p^{p'} + \Vert \bkE (S_{N}   |{\mathcal F}_{0}) \Vert_p^{p'} +   \sum_{k=1}^n
\big \Vert  \sum_{j \geq k+N}  P_0(X_j) \big \Vert_p^{p'} \ \text{ where } \ p'=\min ( 2,p)$.

\end{enumerate}
\end{Proposition}

Let us now consider the following reinforcement of the first part of \eqref{newcondTH22}: there exists a sequence $(u_n)_{n \geq 1}$ of positive reals such that $u_n \geq n$ and
\beq \label{Cond1cobSn}
\sum_{n \geq 2} \frac{ \max_{k \in \{n, u_n \}}\| {\mathbf E}(S_{k} | {\mathcal F}_0) \|^p_p }{n^2 (\log n)^{(t-1)p/2}}< \infty      \,
\, \text{ and } \,
\sum_{n \geq 2}  \frac{ 1 }{n^2 (\log n)^{(t-1)p/2} }  \Big ( \sum_{k=1}^n \Big \| \sum_{j \geq k+u_n}P_0 (X_j) \Big \|^2_p \Big )^{p/2}     < \infty  \, ,
\eeq
and
\beq \label{Cond1cobSn2}
 \sum_{n \geq 2}   \frac{ n^{p/4} }{ n^2 (\log n)^{(t-1)p/2}}
  \Big ( \sum_{k=1}^n \Big \| \sum_{j \geq k+n}P_0 (X_j) \Big \|^2_2 \Big )^{p/4}    < \infty  \, .
  \eeq
With the help of Proposition \ref{approxrnq}, we then obtain the following counterpart to Theorem \ref{interwunous}:

\begin{Theorem} \label{ThNAS} Let $2< p \leq 4$ and $t>2/p$. Assume that $X_0$ belongs to ${\mathbf L}^p$ and that (\ref{Gor}) is satisfied. Assume in addition that the conditions \eqref{Cond1cobSn} and \eqref{Cond1cobSn2} hold and that
\beq \label{condcarre}
\sum_{n \geq 2} \frac{ 1 }{ n^2 (\log n)^{(t-1)p/2} }  \big \Vert  \bkE (S_n^2 | {\mathcal F}_0  ) -\bkE  (S_n^2)
 \big \Vert^{p/2}_{p/2} < \infty \, .\eeq
Then the conclusion of Theorem \ref{interwunous} holds.
\end{Theorem}

In view of applications to  mixingale-like sequences, we give the following results:

\begin{Proposition} \label{directcond} Let $2< p \leq 4$ and $t >2/p$. Assume that $X_0$ belongs to ${\mathbf L}^p$.  If
\beq \label{Cond1cob*}
\sum_{n \geq 2} \frac{n^{p-1}}{n^{2/p } (\log n)^{\frac{(t-1)p}{2 }}} \| {\mathbf E}(X_n| {\mathcal F}_0) \|^p_p < \infty \text{ and }  \sum_{n \geq 2}   \frac{ n^{3p/4}  }{n^2  (\log n)^{\frac{(t-1)p}{2 }}}   \Vert \E(X_n |{\mathcal F}_0 )\Vert^{p/2}_2
 < \infty\, ,
\eeq
then (\ref{Gor}) holds, and the conditions \eqref{Cond1cobSn} and \eqref{Cond1cobSn2} are  satisfied with $u_n=[n^{p/2}]$. In addition $ n^{-1}\E(S_n^2)$
converges to $\sum_{k \in {\mathbf Z}} {\rm Cov } (X_0 , X_k)$
as $n$ tends to infinity.
\end{Proposition}

\begin{Corollary} \label{coralphaphi}
Let $2< p \leq 4$ and $t >2/p$. Assume that $X_0$ belongs to ${\mathbf L}^p$ and that there exists $\gamma\in ]0,1]$ such that
\beq \label{cond1coralphaphi}
 \sum_{n>0}  \frac{ n^{ ( \frac p2 -1 ) ( \frac{1}{\gamma} +1)} }{n^{1/2} (\log n)^{(t-1)p/2}} \| {\mathbf E}(X_n| {\mathcal F}_0) \|^{p/2}_p  < \infty \, ,
\eeq
and
\beq \label{cond2coralphaphi}
 \sum_{n>0}  \frac{n^{(\gamma +1) p/2} }{n^2 (\log n)^{(t-1)p/2}}  \sup_{i \geq j \geq n }\| {\mathbf E}(X_iX_j | {\mathcal F}_0) - {\mathbf E}(X_iX_j)\|_{p/2}^{p/2} < \infty \, .
\eeq
Then the conclusion of Theorem \ref{interwunous} holds with $\sigma^2 = \sum_{k \in \Z} {\rm Cov}(X_0,X_k)$.
\end{Corollary}

\medskip

The next result has a different range of applicability than Corollary \ref{coralphaphi}.

\begin{Corollary} \label{coralpha}
Let $2< p \leq 4$ and $t >2/p$. Assume that $X_0$ belongs to ${\mathbf L}^p$ and that (\ref{Cond1cob*}) holds. Assume in addition that
\beq \label{cond2coralpha}
\sum_{n>0}  \frac{n^{p}}{n^{2} (\log n)^{(t-1)p/2} }  \| X_0 {\mathbf E}(X_n | {\mathcal F}_0)\|^{p/2}_{p/2}  < \infty  \, ,
\eeq
and
\beq \label{cond2coralpha*}
  \sum_{n>0}  \frac{n^{p}}{n^{2} (\log n)^{(t-1)p/2}}  \sup_{i \geq j \geq n }\| {\mathbf E}(X_iX_j | {\mathcal F}_0) - {\mathbf E}(X_iX_j)\|_{p/2}^{p/2} < \infty    \, .
\eeq
Then  the conclusion of Theorem \ref{interwunous} holds with $\sigma^2 = \sum_{k \in \Z} {\rm Cov}(X_0,X_k)$.
\end{Corollary}

\section{Applications} \label{sectionappli}
\setcounter{equation}{0}

\subsection{Application to an example of irreducible Markov chain.} \label{sectinvert}

In this section we apply Theorem \ref{interwunous} to a Markov chain which is a symmetrized version of the Harris recurrent Markov chain defined in Doukhan, Massart and Rio (1994) and that has been considered recently in Rio (2009). Let $E=[-1,1]$ and let $\upsilon$ be a symmetric
atomless law on $E$. The transition probabilities are defined by
\[
Q(x,A)=(1-|x|)\delta_{x}(A)+|x|\upsilon(A) \, ,
\]
where $\delta_{x}$ denotes the Dirac measure at point $x$. Assume that $\theta=\int_{E}%
|x|^{-1}\upsilon(dx)<\infty$. Then there is an unique invariant measure
\beq \label{hypopi}
\pi(dx)=\theta^{-1}|x|^{-1}\upsilon(dx) \, ,
\eeq
and the stationary Markov chain $(\zeta_{i})_i$ is reversible and positively recurrent.

Let $f$ be a measurable function on $E$ and $X_{i}=f(\zeta_{i})$. We denote by
$S_{n}(f)$ the partial sum $S_{n}$. Assuming that the measure $\upsilon$ satisfies
\begin{equation} \label{hyp1upsilon}
\upsilon ([0,t]) \leq c t^{a+1}  \text{ for some $a > (p-2)/2$ and some $c>0$,}
\end{equation}
and that $f$ is an odd function satisfying  $|f(x)|\leq C|x|^{1/2}$ for any $x$ in $E$
with $C$ is a positive constant, Rio (2009) (for $p \in ]2,3]$) and  Merlev\`ede and Peligrad (2012) (for any $p>2$) have shown that
 $\| \max_{1 \leq k \leq n}|S_k(f)| \|_p $ satisfies a
Rosenthal-type inequality. When $p \in ]2,4[$, applying Theorem \ref{ThNAS}, we shall prove that under the same assumptions, $S_n(f)$ satisfies the strong approximation \eqref{strassen} with rate $b_n=n^{1/p}$.
\begin{Corollary}
\label{corMarkovChain} Let $f$ be such that $f(-x)=-f(x)$ for
any $x\in E$. Assume that there exist  $C>0$ and $\gamma \geq 1/2$ such that  $|f(x)|\leq C|x|^{1/2}$ for any $x$ in $E$. Let $p\in]2,4]$ be a real number and assume that (\ref{hyp1upsilon}) holds true. Then $S_n(f)$ satisfies the strong approximation \eqref{strassen} with $\sigma^2 = \sum_{k \in {\mathbf Z}} {\rm Cov } (X_0 , X_k)$ and rate $b_n=n^{1/p} $ if $p \in ]2,4[$ and $b_n=n^{1/4} (\log n)^{3/4 + \varepsilon}$ for any $\varepsilon >0$ if $p=4$.
\end{Corollary}

\begin{Remark}
Let $\beta_{{\bf \zeta}} (n) = 2^{-1}\int_E \Vert Q^n (x, \cdot) - \pi (\cdot ) \Vert \pi (dx)$ where $\Vert\mu(\cdot)\Vert$ denotes the total variation of the signed measure
$\mu$. According to Lemma 2 in Doukhan et al. (1994), the absolute regularity coefficients $\beta_{{\bf \zeta}} (n)$ of the sequence $(\zeta_i)_i$ are exactly of order $n^{-a}$. Therefore, for $p\in ]2,4[$, as soon as $\gamma$ is big enough, $S_n(f)$ satisfies the almost sure invariance principle with the rate $n^{1/p}$ even if the absolute regularity coefficients of the Markov chain $(\zeta_{i})_i$ do not satisfy $\sum_{n \geq 1} \beta_{{\bf \zeta}} (n)< \infty$ which corresponds to the ergodicity of degree two (see Nummelin (1984)).

Notice also that an application of Theorem 2.1 in Merlev\`ede and Rio (2012) would require $a > p-1$ if $ p \in ]2,3]$ to get the rate $n^{1/p}$ (up to some logarithmic terms) in the almost sure invariance principle for $S_n(f)$.
\end{Remark}

\subsection{Symmetric random walk on the circle} \label{sectcircle}
 Let $K$ be the Markov kernel defined by $$Kf (x) =
\frac{1}{2}( f (x+a) + f (x-a) )$$ on the torus $\R/\Z$, with $a$
irrational in $[0,1]$. The Lebesgue-Haar  measure $m$  is the unique
probability which is invariant by $K$. Let $(\xi_i)_{i\in \Z}$ be
the stationary Markov chain with transition kernel $K$  and
invariant distribution $m$. For $f \in {\bf L}^2(m)$, let
\begin{equation}\label{defSnf}
X_k=f(\xi_k)-m(f) \, .
\end{equation}
This example has been considered by Derriennic and Lin (2001) who showed (see their section 2) that the central
limit theorem holds for $n^{-1/2} \sum_{k=1}^nX_k$ as soon as the series
of covariances
\begin{equation}\label{sericov}
\sigma^2(f) = m((f-m(f))^2) + 2 \sum_{n>0} m( f K^n(f-m(f)))
\end{equation}
is convergent, and that the limiting distribution is ${\mathcal
N}(0,\sigma^2(f))$. In fact the convergence of the series in
(\ref{sericov}) is equivalent to
\begin{equation}\label{Paroux}
 \sum_{k\in{\Z}^*}   \frac{|\hat f (k)|^2}{d (ka, {\Z})^2} < \infty
 \, ,
\end{equation}
where $\hat f(k)$ are the Fourier coefficients of $f$
and $d (ka, {\Z})=\min_{i \in \Z} |ka-i|$. Hence, for
any irrational number $a$, the criterion (\ref{Paroux}) gives a
class of function $f$ satisfying the central limit theorem, which
depends on the sequence $((d (ka, {\Z}))_{k \in{\Z}^*}$. Note that a
function $f$ such that
\begin{equation}\label{debilmental}
\liminf_{k \rightarrow \infty} k |\hat f(k)|>0 \, ,
\end{equation}
does not satisfies (\ref{Paroux}) for any irrational number $a$.
Indeed, it is well known from the theory of continued fraction that
if $p_n/q_n$ is the $n$-th convergent of $a$, then $|p_n-q_n
a|<q_n^{-1}$, so that $d(ka, \Z)<k^{-1}$ for an infinite number of
positive integers $k$. Hence, if (\ref{debilmental}) holds, then
$|\hat f(k)|/d(ka, \Z)$ does not even tend to zero as $k$ tends to
infinity.

Our aim in this section is to give conditions on $f$ and on the
properties of the irrational number $a$ ensuring rates of convergence in the almost sure sure invariance principle. Let us then introduce the following definitions: $a$ is said to be \textit{badly approximable}
by rationals if
\beq \label{badly} d(ka, \Z) \geq c(a) |k|^{-1} \ \text{ for some positive constant $c(a)$}.
\eeq
An irrational number is badly approximable iff the terms $a_n$ of its continued fraction are bounded. In particular, the quadratic irrationals are badly approximable. However, note that the set of numbers in $[0,1]$ satisfying \eqref{badly} has Lebesgue measure $0$. A much bigger set is the following: $a$ is said to be \textit{badly approximable in the weak sense} by rationals if for any positive $\varepsilon$,
\beq\label{badlyweak} \text{ the inequality $d(ka, \Z) < |k|^{-1-\varepsilon}$ has only finitely many
solutions for $k \in \Z$.}
\eeq

\noindent From Roth's theorem the algebraic numbers are badly
approximable in the weak sense (cf. Schmidt (1980)). Note also that the set of  numbers in $[0,1]$ satisfying \eqref{badlyweak}, has Lebesgue measure $1$.
Let us note that in Section 5.3 of Dedecker and Rio (2008), it is proved that the
condition (\ref{Paroux}) (and hence the central limit theorem for
$n^{-1/2} \sum_{k=1}^nX_k$) holds for any number $a$ satisfying \eqref{badlyweak}
as soon as
\begin{equation}\label{antideb}
\sup_{k\not= 0} |k|^{1+\varepsilon} |\hat f (k)| < \infty \quad
\text{for some positive $\varepsilon$.}
\end{equation}
Note that, in view of (\ref{debilmental}), one cannot take
$\varepsilon=0$ in  the condition (\ref{antideb}).

\begin{Corollary}\label{circle0}
Let $X_k$ be defined by \eqref{defSnf}. Suppose that $a$ satisfies \eqref{badly}. Let $p \in ]2,4]$ and assume that for some positive $\varepsilon$,
\begin{equation} \label{condirra0}
\sup_{k\not= 0} |k|^{ s} \big ( \log ( 1 + |k|) \big )^{1+\varepsilon} |\hat f (k)| < \infty \quad
\text{ where }s=\frac{\sqrt{1+4p(p-2)}}{p} - \frac{3}{p} + 2  .
\end{equation}
Then $S_n(f)$ satisfies the strong approximation \eqref{strassen} with $\sigma^2 = \sum_{k \in {\mathbf Z}} {\rm Cov } (X_0 , X_k)$ and rate $b_n=n^{1/p}  \log n  $.
\end{Corollary}

When the condition on $a$ is weaker, we obtain:

\begin{Corollary}\label{circle}
Let $X_k$ be defined by \eqref{defSnf}. Suppose that $a$ satisfies \eqref{badlyweak}. Let $p \in ]2,4]$ and assume that for $s$ defined in \eqref{condirra0} and some positive $\varepsilon$,
\begin{equation*} \label{condirra}
\sup_{k\not= 0} |k|^{ s +\varepsilon} |\hat f (k)| < \infty \, .
\end{equation*}
Then $S_n(f)$ satisfies the strong approximation \eqref{strassen} with $\sigma^2 = \sum_{k \in {\mathbf Z}} {\rm Cov } (X_0 , X_k)$ and rate $b_n=n^{1/p} $ if $p \in ]2,4[$ and $b_n=n^{1/4} (\log n)^{3/4 + \delta}$ for any $\delta >0$ if $p=4$.
\end{Corollary}
\begin{Remark}\label{rkrw} Applying Corollary \ref{circle} with $p$ close enough to 2, we derive that if the function $f$ satisfies \eqref{antideb} then,
enlarging $\Omega$ if
necessary, there exists a sequence $(Z_i)_{i
\geq 1}$ of iid gaussian random variables with zero mean and
variance $\sigma^2$ such that, for some $\eta >0$,
\beq \label{rem***}
\sum_{i=1}^n (X_i -Z_i)  = o (n^{1/2 - \eta}  ) \text{ almost surely, as $n\rightarrow \infty$},
\eeq
which could be also deduced from Theorem 1 in Eberlein (1986). As a consequence of \eqref{rem***}, the  weak invariance principle as well as the almost sure invariance principle hold true under \eqref{antideb}. Note also that if $f $ is three times differentiable then $\sum_{i=1}^n X_i$ satisfies the strong approximation \eqref{strassen} with rate $b_n = n^{1/4} ( \log n)^{3/4 + \delta}$ for any $\delta >0$.
\end{Remark}

\subsection{Application to a class of weak dependent sequences} \label{secttau}
In this section  we give rates of
convergence in the almost sure invariance principle for a stationary
sequence $(X_i)_{i \in {\mathbf Z}}$ satisfying some weak dependence conditions  specified below.

\begin{Definition} Let $\Lambda_1({\mathbf R})$ be the set of the  functions $f$ from ${\mathbf R}$ to ${\mathbf R}$ such that $|f(x)-f(y)|\leq |x-y|$.
For any $\sigma$-algebra
${\mathcal F }$ of ${\mathcal A}$ and any real-valued integrable random
variable $X$, we consider the coefficient $\theta ({\mathcal F}, X )$ defined by
 \beq \label{deftheta1} \theta ({ \mathcal F} , X ) =
\sup_{f \in \Lambda_1(\R)}\| \E ( f(X) | { \mathcal F} ) -\E ( f(X))
\|_1    \, .\eeq
\end{Definition}
We now define the coefficients $\gamma(n)$, $\theta_2(n)$ and $\lambda_2(n)$ of the sequence
$(X_i)_{i \in {\mathbf Z}}$.
\begin{Definition}  For any positive integer $k$, define \beq
\label{deftheta2} \theta_2 (n) = \sup_{i \geq j \geq n}
\max \{ \theta ({ \cal F}_0 , X_i + X_j ) , \theta ({ \cal F}_0 , X_i - X_j ) \} \quad  \text{and} \quad
 \gamma (n) =  \| \E ( X_n| { \cal F}_0 ) \|_{1}\, . \eeq
 Let now \beq \label{deftilde}
\lambda_2 (n) = \max ( \theta_2 (n) ,  \gamma (n)) \, .\eeq
\end{Definition}

\begin{Definition} \label{defquant}
 For any integrable random variable
$X$, define the ``upper tail'' quantile function $Q_X $ by $
Q_X (u) = \inf \left \{  t \geq 0 : {\mathbf P} \left(|X| >t \right) \leq
u\right \} $. Note that, on the set $[0, {\mathbf P}(|X|>0)]$, the function
$H_X: x \rightarrow \int_0^x Q_X(u) du$ is an absolutely continuous
and increasing function with values in $[0, {\mathbf E}|X|]$. Denote
by  $G_X$ the inverse of $H_X$.
\end{Definition}

\begin{Corollary}
\label{cortheta}
Let $2< p \leq 4$ and $t >2/p$.  Assume that $X_0$ belongs to ${\mathbf L}^p$. Let  $Q=
Q_{X_0}$, and $G= G_{X_0}$.  Assume in addition that
 \beq \label{Condtheta}
\sum_{n \geq 2} \frac{n^{p-1}}{n^{2/p } (\log n)^{\frac{(t-1)p}{2 }}}  \int_0^{\lambda_2(n)} Q^{p-1} \circ G (u) du
 < \infty \, .\eeq
Then the conclusion of Theorem \ref{interwunous} holds  with $\sigma^2 = \sum_{k \in {\mathbf Z}} {\rm Cov} (X_0,X_k)$.
\end{Corollary}
Denote by ${\mathcal F}_{\ell} = \sigma ( X_i, i \leq \ell)$ and by
and ${\mathcal G}_k=\sigma ( X_i, i \geq k)$. Notice that if $\alpha$ denote the usual strong mixing coefficient of Rosenblatt (1956) of $X$ defined by
\begin{equation*}\label{defalpharosen}
 \alpha(n)= \sup_{\ell \in {\mathbf Z}}\alpha({\mathcal F}_{\ell},{\mathcal G}_{n + \ell}) \text{ for $n\geq 0$} \, ,
\end{equation*}
where $  \alpha({\mathcal F}, {\mathcal G})= \sup_{A \in {\mathcal F}, B \in {\mathcal G}}|{\mathbf P}(A \cap B)-{\mathbf P}(A){\mathbf P}(B)| $, then  according to Lemma 1 in Dedecker and Doukhan (2003), condition \eqref{Condtheta} is implied by
$$
\sum_{n \geq 2} \frac{n^{p-1}}{n^{2/p } (\log n)^{\frac{(t-1)p}{2 }}}  \int_0^{\alpha(n)} Q^{p} (u) du
 < \infty \, .$$
 As a consequence, it follows that, if for $r>p$,
 \beq \label{condalphafort}
\sup_{x >0} x^r {\mathbf P}(|X_0| > x) < \infty \quad \text{ and } \quad \sum_{n \geq 2} \frac{n^{p-1}}{n^{2/p  } }  (\alpha(n) )^{(r-p)/r} < \infty\, ,
\eeq
then (\ref{strassen}) holds with $b_n = n^{1/p} (\log n)$. Therefore, Corollary \ref{coralphafaible} improves on Shao and Lu's result (1987), which requires: $\sum_{n=1}^\infty
 (\alpha(n) )^{(r-p)/(rp)} < \infty$. Notice however that \eqref{condalphafort}  is stronger than $\sum_{n=1}^\infty
n^{p-2} (\alpha(n) )^{(r-p)/r} < \infty$ which is the condition obtained by Merlev\`ede and Rio (2012) to get the rate $n^{1/p}$ (up to  
logarithmic terms) in (\ref{strassen}), but only for $p \in ]2,3]$. Therefore, compared to their Theorem 2.1, Corollary \ref{cortheta} allows  better  rates than $n^{1/3}$.

\subsubsection{Application to $\tau$-dependent sequences}
As we shall see Corollary \ref{cortheta} is well adapted to obtain rates of
convergence in the almost sure invariance principle for functions of $\tau$-dependent sequences. Before stating the result, some definitions are needed.

\begin{Definition}
Let $\Lambda_1({\mathbf R})$ be the set of the  functions $f$ from ${\mathbf R}$ to ${\mathbf R}$ such that $|f(x)-f(y)|\leq |x-y|$.
Let $\Lambda_1({\mathbf R}^2)$ be the set of functions $f$
from ${\mathbf R}^2$ to ${\mathbf R}$ such that
$$
  |f(x_1, y_1)-f(x_2, y_2)|\leq \frac 12 |x_1-y_1|+ \frac 12 |x_2-y_2| \, .
$$
Define the dependence coefficients $\tau_{1,{\bf Y}}$ and $\tau_{2,{\bf Y}}$ of the sequence $(Y_i)_{i \in {\mathbf Z}}$ by
\begin{eqnarray*}
\tau_{1,{\bf Y}} (k) &=& \Big \|\sup_{f \in \Lambda_1({\mathbf R})}\Big |{\mathbf E}(f(Y_k)|{\mathcal F}_0)-{\mathbf E}(f(Y_k))\Big | \Big \|_1 \, , \\
 \tau_{2,{\bf Y}} (k) &=& \max \Big \{ \tau_{1,{\bf Y}}(k),
  \sup_{i > j \geq k} \Big \| \sup_{f \in \Lambda_1({\mathbf R}^2)}\Big |{\mathbf E}(f(Y_i, Y_j)|{\mathcal F}_0)-{\mathbf E}(f(Y_i, Y_j))\Big | \Big \|_1 \Big \}\, .
\end{eqnarray*}
\end{Definition}
Many examples of $\tau$-dependent sequences are given in Dedecker and Prieur (2005).

\medskip

We now define the classes of functions which are adapted to this kind of dependence.


\begin{Definition}
Let $c$ be any concave function from ${\mathbf R}^+$ to ${\mathbf R}^+$, with $c(0)=0$. Let ${\mathcal L}_c$
be the set of functions $f$ from ${\mathbf R}$ to ${\mathbf R}$ such that
$$
 |f(x)-f(y)| \leq K c(|x-y|), \quad \text{for some positive $K$.}
$$
\end{Definition}

An application of Corollary \ref{cortheta} gives:

\begin{Corollary} \label{cortau}
Let $f \in {\mathcal L}_c$, and
let $X_k= f(Y_k) - {\mathbf E}(f(Y_k))$. Let $2< p \leq 4$ and $t >2/p$. Assume that $X_0$ belongs to ${\mathbf L}^p$. Let  $Q=
Q_{X_0}$, and $G= G_{X_0}$.  Assume in addition that
 \beq \label{Condtau}
\sum_{n \geq 2} \frac{n^{p-1}}{n^{2/p  } (\log n)^{\frac{(t-1)p}{2 }}} \int_0^{c(\tau_{2,{\bf Y}}(n))} Q^{p-1} \circ G (u) du
 < \infty \, .\eeq
Then the conclusion of Theorem \ref{interwunous} holds  with $\sigma^2 = \sum_{k \in {\mathbf Z}} {\rm Cov} (X_0,X_k)$.
\end{Corollary}
Notice that if  for $r>p$ and $t >2/p$, $$
\sup_{x >0} x^r {\mathbf P}(|X_0| > x) < \infty \quad \text{ and } \quad \sum_{n \geq 2} \frac{n^{p-1}}{n^{2/p } (\log n)^{\frac{(t-1)p}{2 }}}(c(\tau_{2,{\bf Y}}(k)) )^{(r-p)/(r-1)} < \infty\, ,
$$
the condition \eqref{Condtau} is satisfied.

\medskip

\noindent{\it Proof of Corollary \ref{cortau}.} Let us first prove that the condition (\ref{Condtau}) implies the condition (\ref{Condtheta}). Let $X_k=f(Y_k) - {\mathbf E}(f(Y_k))$. Applying the coupling result given in  Dedecker and Prieur (2005, Section 7.1) (see also Proposition 4 in R\"uschendorf (1985)),
we infer that there exists $\bar Y_n$ distributed as $Y_n$ and independent
of ${\mathcal F}_0$ such that
$$
{\mathbf E}(|Y_n-\bar Y_n|)= \tau_{1,{\bf Y}}(n) \leq \tau_{2,{\bf Y}}(n)\, .
$$
In the same way, for $n \leq i<j$ there exists $(Y_i^*, Y_j^*)$ distributed as $(Y_i, Y_j)$ and independent of ${\mathcal F}_0$
such that
$$
  \frac12{\mathbf E}(|Y_i-Y_i^*|+ |Y_j-Y_j^*|)=
  \Big \| \sup_{h \in \Lambda_1({\mathbf R}^2)}\Big |{\mathbf E}(h(Y_i, Y_j)|{\mathcal F}_0)-{\mathbf E}(h(Y_i, Y_j))\Big | \Big \|_1 \leq \tau_{2,{\bf Y}}(i) \leq \tau_{2,{\bf Y}}(n)\, .
$$

Clearly
$$
 \gamma(n) =\|{\mathbf E}(f(Y_n)|{\mathcal F}_0)-{\mathbf E}(f(Y_n))\|_{1} \leq \|f(Y_n)-f(\bar Y_n)\|_{1}\, .
$$
Consequently, if $f \in {\mathcal L}_c$, one has
$$
  \gamma(n) \leq K{\mathbf E}(c(|Y_n-\bar Y_n|))\leq Kc(\|Y_n-\bar Y_n\|_1)=Kc(\tau_{1,{\bf Y}}(n)) \, .
$$

In the same way, if $g$ is in $\Lambda_1({\mathbf R})$,
$$
 \|{\mathbf E}(g(X_i+X_j)|{\mathcal F}_0)-{\mathbf E}(g(X_i+X_j))\|_1 \leq {\mathbf E}(|f(Y_i)-f(Y_i^*)|+ |f(Y_j)-f(Y_j^*)|)\, .
$$
Hence, if  $f\in {\mathcal L}_c$,
$$
\|{\mathbf E}(g(X_i+X_j)|{\mathcal F}_0)-{\mathbf E}(g(X_i+X_j))\|_1 \leq 2 K c(\tau_{2,{\bf Y}}(n)) \, .
$$
Note that the same inequalities hold with $X_i-X_j$ instead of $X_i+X_j$.

As a consequence, we obtain that if $f \in {\mathcal L}_c$, then $\lambda_2(n) \leq 2K c(\tau_{2,{\bf Y}}(n))$. Hence, Corollary \ref{cortau}  follows from Corollary \ref{cortheta}.  $\diamond$

\medskip

\noindent{\it Example: Autoregressive Lipschitz model.} Let us give an example of an iterative Lipschitz model,
which may fail to be irreducible and to which Corollary \ref{cortau} applies. For $\delta$ in
$[0,1[$ and $C$ in $]0,1]$, let ${\mathcal L} (C,\delta)$ be the
class of 1-Lipschitz functions $h$ which satisfy
$$
h (0) = 0 \ \mbox{ and }\ |h' (t)| \leq 1-C ( 1+ |t| )^{-\delta} \
\mbox{ almost everywhere.}
$$
Let $(\varepsilon_i)_{i\in { \mathbb {Z} } }$ be a sequence of
i.i.d. real-valued random variables. For  $S \geq 1$, let $ARL(C,\delta,S)$ be the class of Markov
chains on ${\mathbf R}$ defined by
\begin{equation}
Y_n=h(Y_{n-1})+\varepsilon_n \ \mbox{ with }\ h \in {\mathcal L}
(C,\delta)
 \ \mbox{ and }\ {\mathbf E}(|\varepsilon_0|^S)< \infty \, .
\end{equation}
For this model, there exists an unique invariant probability measure $\mu$ such that $\mu (|x|^{S-\delta}) < \infty$ (see Proposition 2 of
Dedecker and Rio (2000)). In addition starting from the inequality (7.7) in Dedecker and Prieur (2005) and arguing as in Dedecker and Rio (2000), one can prove that $\tau_{2,{\bf Y}}(n)= O (n^{(\delta +1 -S)/\delta}) $ if $S>1+\delta$. Therefore an application of Corollary \ref{cortau} leads directly to the following result:
\begin{Corollary} \label{corarl}
Assume that $(Y_i)_{ i\in \mathbf{Z} } $ belongs to
$ARL(C,\delta,S)$. Let $f$ be some H\"older function of order $\gamma \in ]0,1]$,  that is  $|f(x) - f(y)| \leq K|x-y|^{\gamma}$ for some $K>0$. Let $X_i = f (Y_i)- {\mathbf{E}} (f(Y_i))$ and $S_n(f)=\sum_{k=1}^n X_k$. If for some $p\in ]2,4]$,
\begin{equation} \label{condfap}
S>1+\delta \, \, \text{ and } \, \, \frac{(S-1-\delta) ( S-\delta - \gamma p)}{S - \gamma - \delta} > \frac{\delta}{\gamma} \Big ( p - \frac 2p \Big )\, ,\end{equation}
then $S_n(f)$ satisfies the strong approximation \eqref{strassen} with $\sigma^2 = \sum_{k \in {\mathbf Z}} {\rm Cov } (X_0 , X_k)$ and rate $b_n=n^{1/p} $ if $p \in ]2,4[$ and $b_n=n^{1/4} (\log n)^{3/4 + \varepsilon}$ for any $\varepsilon >0$ if $p=4$.
\end{Corollary}
Notice that the condition $S >p+  \delta \big (1 + \gamma^{-1} (p -2/p) \big ) $ implies the condition (\ref{condfap}) (both conditions are identical if $\gamma=1$, that is when $g$ is Lipschitz).

An element of  $ARL(C,\delta,\eta)$ may fail to be irreducible and then strongly mixing in
the general case. However, if the common distribution of the
$\varepsilon_i$'s has an absolutely continuous component which is
bounded away from $0$ in a neighborhood of the origin, then the
chain is irreducible and fits in the example of Tuominen and Tweedie
(1994), Section 5.2. In this case, the rate of ergodicity can be
derived from Theorem 2.1 in Tuominen and Tweedie (1994).

\subsection{Application to $\alpha$-dependent sequences} \label{sectionappliwms} In this section  we want to consider a weaker coefficient than the Rosenblatt strong mixing coefficient defined by \eqref{defalpharosen}, and which may computed for instance for many Markov chains associated to dynamical systems that fail to be strongly mixing. We start with the definition of the $\alpha$-dependent coefficients.

\begin{Definition}
For any integrable random variable $X$, let us write
$X^{(0)}=X- \bkE(X)$.
For any random variable $Y=(Y_1, \cdots, Y_k)$ with values in
${\mathbf R}^k$ and any $\sigma$-algebra $\F$, let
\[
\alpha(\F, Y)= \sup_{(x_1, \ldots , x_k) \in {\mathbf R}^k}
\left \| \bkE \Big(\prod_{j=1}^k (\I_{Y_j \leq x_j})^{(0)} \Big | \F \Big)^{(0)} \right\|_1.
\]
For the sequence ${\bf Y}=(Y_i)_{i \in {\mathbf Z}}$, let \begin{equation}
\label{defalpha} \alpha_{k, {\bf Y}}(0) =1 \text{ and }\alpha_{k, {\bf Y}}(n) = \max_{1 \leq l \leq
k} \ \sup_{ n\leq i_1\leq \ldots \leq i_l} \alpha(\F_0,
(Y_{i_1}, \ldots, Y_{i_l})) \text{ for $n>0$}.
\end{equation}Let $\Lambda_1({\mathbf R})$ be the set of the  functions $f$ from ${\mathbf R}$ to ${\mathbf R}$ such that $|f(x)-f(y)|\leq |x-y|$.
\end{Definition}
Notice that $\alpha_{k, {\bf Y}}(n) \leq
\alpha (n)$ for any positive $n$, where $\alpha(n)$ is the strong mixing coefficient of Rosenblatt of ${\bf Y}$ as defined by \eqref{defalpharosen}. For examples of Markov chains satisfying  $\lim_{n \rightarrow \infty}\alpha_{k, {\bf Y}}(n) = 0 $  and which are not strongly mixing in the sense of Rosenblatt, see the section \ref{sectweakMC}.

\medskip

We now define the classes of functions which are adapted to  this kind of weak dependence.

\begin{Definition}
\label{defclosedenv} Let $\mu$ be the probability distribution of a
random variable $X$. If $Q$ is an integrable quantile function (see definition \ref{defquant}),
let $\tMon( Q, \mu)$ be the set of functions $g$ which are
monotonic on some open interval of ${\mathbf R}$ and null
elsewhere and such that $Q_{|g(X)|} \leq Q$. Let $\tMonm( Q,
\mu)$ be the closure in ${\mathbf L}^1(\mu)$ of the set of
functions which can be written as $\sum_{\ell=1}^{L} a_\ell
f_\ell$, where $\sum_{\ell=1}^{L} |a_\ell| \leq 1$ and $f_\ell$
belongs to $\tMon( Q, \mu)$.
\end{Definition}

For functions of $\alpha$-dependent sequences, the following result holds:

\begin{Corollary} \label{coralphafaible}
Let $2< p \leq 4$ and $t >2/p$. Let  $X_i = f(Y_i) - \bkE ( f(Y_i))$ where $f$ belongs to $\tMonm( Q, P_{Y_0})$ (here,
$P_{Y_0}$ denotes the distribution of $Y_0$) with $Q^p$ integrable. Let
 $\alpha_{2, {\bf Y}}(n) $ be defined as in
(\ref{defalpha}).  Assume that
 \beq \label{Condstrong}
\sum_{n \geq 2} \frac{n^{p-1}}{n^{2/p} (\log n)^{\frac{(t-1)p}{2 }}} \int_0^{\alpha_{2, {\bf Y}}(n)} Q^p(u) du  < \infty \, ,\eeq
Then the conclusion of Theorem \ref{interwunous} holds  with $\sigma^2 = \sum_{k \in {\mathbf Z}} {\rm Cov} (X_0,X_k)$.
\end{Corollary}

When $p$ is close to 2, the condition  (\ref{Condstrong}) is close to the condition: $
\sum_{k\geq
1} \int_0^{\alpha_{2, {\bf Y}}(k)} Q^2(u) du
 < \infty $
which is the best known condition (and optimal in a sense) for the strong invariance principle of $\alpha$-dependent sequences (see Theorem 1.13 of Dedecker, Gou\"ezel and Merlev\`ede (2010)). However when $p \in ]2,3]$, Theorem 2.1 in Merlev\`ede and Rio (2012) provides a sharper condition than \eqref{Condstrong}. As a counterpart, our Corollary \ref{coralphafaible} allows rates of convergence of order $n^{1/p}$ (up to logarithmic terms) with $p \in ]3,4]$ in the almost sure invariance principle that are not reached in  Merlev\`ede and Rio's paper.

\subsubsection{Application to  functions of  Markov chains associated to  intermittent maps} \label{sectweakMC}

For  $\gamma$ in $]0, 1[$, we consider the intermittent map
$T_\gamma$ from $[0, 1]$ to $[0, 1]$, which is a modification of the
Pomeau-Manneville map (1980):
$$
   T_\gamma(x)=
  \begin{cases}
  x(1+ 2^\gamma x^\gamma) \quad  \text{ if $x \in [0, 1/2[$}\\
  2x-1 \quad \quad \quad \ \  \text{if $x \in [1/2, 1]$} \, .
  \end{cases}
$$
We denote by $\nu_\gamma$ the unique $T_\gamma$-invariant
probability measure on $[0, 1]$ which is absolutely continuous with
respect to the Lebesgue measure. We denote by $K_\gamma$ the
Perron-Frobenius operator of $T_\gamma$ with respect to
$\nu_\gamma$. Recall that for any bounded measurable functions $f$ and $ g$,
$$
\nu_\gamma(f \cdot  g\circ T_\gamma)=\nu_\gamma(K_\gamma(f) g) \, .
$$
Let $(Y_i)_{i \geq 0}$ be a stationary Markov chain with invariant
measure $\nu_\gamma$ and transition Kernel $K_\gamma$.  Applying Corollary \ref{coralphafaible}, we shall see that for $f$ belonging to a certain class of functions defined  below, $\sum_{k=1}^n (f(Y_i)-\nu_\gamma(f))$ satisfies the strong approximation principle \eqref{strassen} with rate $b_n= n^{1/p} ( \log n)$.
\begin{Definition}
A function $H$ from ${\mathbb R}_+$ to $[0, 1]$ is a tail
function if it is non-increasing, right continuous, converges
to zero at infinity, and $x\rightarrow x H(x)$ is integrable.
If $\mu$ is a probability measure on $\mathbb R$ and $H$ is a
tail function, let $\Mon(H, \mu)$ denote the set of functions
$f:\R\to \R$ which are monotonic on some open interval and null
elsewhere and such that $\mu(|f|>t)\leq H(t)$. Let $\Monm(H,
\mu)$ be the closure in ${\bf L}^1(\mu)$ of the set of
functions which can be written as $\sum_{\ell=1}^L a_\ell
f_\ell$, where $\sum_{\ell=1}^L |a_\ell| \leq 1$ and $f_\ell\in
\Mon(H, \mu)$.
\end{Definition}
\begin{Corollary} \label{ASmapB}  Let $(Y_i)_{i \geq 1}$ be a stationary Markov chain with
transition kernel $K_{\gamma}$ and invariant measure  $\nu_\gamma$. Let $p \in ]2,4]$ and let $H$ be a tail function with
\begin{equation}\label{lilcond}
\int_0^{\infty} x^{p-1} (H(x))^{\frac{1-\gamma  \delta }{1-\gamma}} dx <\infty \text{ where $\delta= p+1 -2/p$}\,.
\end{equation}
Then, for any $f \in \Monm(H, \nu_{\gamma})$, the series
\[\sigma^2= \nu_\gamma((f-\nu_\gamma(f))^2)+ 2
\sum_{k>0} \nu_\gamma ((f-\nu_\gamma(f))f\circ T_\gamma^k)
\]
converges absolutely to some nonnegative number, and $\sum_{k=1}^n (f(Y_i)-\nu_\gamma(f))$ satisfies the strong invariance principle \eqref{strassen} with $b_n= n^{1/p} ( \log n)$.
\end{Corollary}
To prove this  corollary, it suffices to see that (\ref{lilcond}) implies 
(\ref{Condstrong}) with $t=1$. In this purpose,
 we use Proposition 1.17 in Dedecker, Gou\"{e}zel and
Merlev\`{e}de (2010) stating that there exist two  positive constant $B, C$ such
that, for any $n>0$, $B n^{(\gamma -1)/\gamma} \leq \alpha_{2,{\bf Y}}(n) \leq C n^{(\gamma -1)/\gamma}$, together with their computations page 817.

Note that Corollary \ref{ASmapB} can  be extended to functions of  Markov chains associated to generalized Pomeau-Manneville
maps (or GPM maps) of parameter $\gamma \in (0,1)$ as defined in Dedecker, Gou\"{e}zel and
Merlev\`{e}de (2010). Notice also that when $f$ is a bounded variation function, Corollary \ref{ASmapB} applies as soon as $\gamma \leq \delta^{-1}$. Therefore when $\gamma < 3/10$, we obtain better rates than the one obtained by Merlev\`ede and Rio (2012, Corollary 3.1). In particular, if $\gamma\leq 2/9$, we obtain the rate
$b_n=n^{1/4} \log n$ in \eqref{strassen}.

\section{Proofs} \label{sectproofs}
\setcounter{equation}{0}
Recall that $d_0 = \sum_{j \geq 0} P_0 (X_j)$, so it  is an element of $H_0 \ominus H_{-1}$ and by (\ref{Gor}), it is in ${\mathbf L}^p$. Recall also that $M_n = \sum_{i=1}^n d_0 \circ T^i$ and let $R_n = S_n - M_n$.
\subsection{Proof of Theorem \ref{interwunous}}
We start the proof by stating the following proposition concerning the almost sure convergence of $R_n$ (its proof will be given later).
\begin{Proposition} \label{propproofth22}
Let $p >1$. Assume that $X_0$ belongs to ${\mathbf L}^p$ and that (\ref{Gor}) is satisfied. Let $(\psi(n))_{n \geq 1}$ be a positive and nondecreasing sequence such that there exists a positive constant $C$ satisfying $\psi(2n) \leq C \psi(n) $ for all $n \geq 1$. Assume that \beq \label{psintercond}
\sum_{n \geq 2}\frac{\| R_{n}\|_{q}^q }{n ( \psi(n))^{q}} < \infty  \, \text{ for some $q \in [1,p]$ and } \, \sum_{n \geq 2} \frac{1}{ ( \psi (n) )^{p}}\bigg(\sum_{k=1}^{n} \frac{\|\E (S_{k} |{ \mathcal F}_0)\|_{p}}{k^{1+1/p}}\bigg)^{p} < \infty\, ,
\eeq
then $ R_n = o (\psi(n) )$  almost surely.
\end{Proposition}
With the help of the above proposition, we prove now that under the first part of \eqref{newcondTH22},
\beq \label{psinter}
R_n = o \big ( n^{1/p} ( \log n)^{(t+1)/2}  \big ) \text{ almost surely.}
\eeq
Taking $\psi(n)=  n^{1/p} ( \log n)^{(t+1)/2}$,  we observe that $(\psi(n))_{n \geq 1}$ satisfies the assumptions of Proposition \eqref{propproofth22}. With the above selection of $\psi(n)$,  we infer that a suitable application of H\"older's inequality implies that the second part of \eqref{psintercond} holds provided that the first part of condition \eqref{newcondTH22} is satisfied. 

We prove now that the first part of condition \eqref{newcondTH22} implies the first part of \eqref{psintercond} with $\psi(n)=  n^{1/p} ( \log n)^{(t+1)/2}$ and $q=2$. 
Note that the first part of \eqref{newcondTH22} implies that 
\begin{equation}\label{PelUt}
\sum_{k>0} \frac{\Vert \E (S_{k} |{ \mathcal F}_0)\Vert_2 }{k^{3/2}} < \infty \, .
\end{equation}
Applying  Proposition 1 in Merlev\`ede {\it et al.} (2012), we derive that:
$$
\| R_{n}\|_{2} \ll n^{1/2} \sum_{k \geq n} \frac{\Vert \E (S_{k} |{ \mathcal F}_0)\Vert_2 }{k^{3/2}} 
$$
(notice that under \eqref{Gor}, the approximating martingale considered in the paper by Merlev\`ede {\it et al.} is almost surely equals to $\sum_{k=1}^n d_0 \circ T^i$ where $d_0 = \sum_{j \geq 0} P_0 (X_j)$). Next, using H\"older's inequality, we derive that for any $\gamma \in ]0,1-2/p[$
\begin{multline*}
\sum_{n \geq 2} \frac{1}{ n^{1+2/p} ( \log n)^{t+1}}\| R_{n}\|_{2}^2  \ll \sum_{n \geq 2}\frac{1}{n^{2/p} ( \log n)^{t+1}}\bigg(\sum_{k \geq n} \frac{\Vert \E (S_{k} |{ \mathcal F}_0)\Vert_2 }{k^{3/2}} \bigg)^2 \\
  \ll \sum_{n \geq 2}\frac{1}{  n^{\gamma +2/p} ( \log n)^{t+1}}\sum_{k \geq n} \frac{\Vert \E (S_{k} |{ \mathcal F}_0)\Vert_2^2 }{k^{2-\gamma}}  \ll \sum_{k \geq 2} \frac{\Vert \E (S_{k} |{ \mathcal F}_0)\Vert_2^2 }{ n^{1+2/p} ( \log n)^{t+1}} \, ,
\end{multline*}
which is finite under the first part of condition \eqref{newcondTH22} (to see this it suffices to apply H\"older's inequality).

Therefore, due to the almost sure convergence \eqref{psinter}, to complete the proof of Theorem \ref{interwunous}, it suffices to notice that under the second part of condition \eqref{newcondTH22}, $(M_n)_{n \geq 1}$ satisfies the condition (\ref{condcarremart}) of  Proposition \ref{propmart}
 given in Appendix with $\psi(n)=n^{2/p} ( \log n)^{t} $. Hence, enlarging $\Omega$ is necessary, there exists a sequence $(Z_i)_{i \geq 1}$ of iid centered  Gaussian variables with variance $ {\mathbf E}(d_0^2)$ such that \eqref{strassen} holds with $b_n = n^{1/p} ( \log n)^{(t+1)/2} $. In addition, note that \eqref{PelUt} is a 
 sufficient  condition for $n^{-1} \bkE (S_n^2)$ to converge (see for instance Theorem 1 in Peligrad and Utev (2005)).

 \medskip

To end the proof of the theorem, it remains to prove Proposition \ref{propproofth22}. With this aim, we first notice that due to the properties of monotonicity of the sequence $(\psi(n))_{n \geq 1}$, the almost sure convergence \eqref{psinter} will follow if we can prove that for any $\lambda >0$,
\begin{equation} \label{but1thminter}
\sum_{r \geq 1}{\mathbf P}\Big (\max_{1\leq i \le 2^r} |S_i-M_i|  \geq  \lambda \psi(2^r) \Big ) < \infty \, .
\end{equation}
Let $q\in [1,p]$. Applying inequality (8) of Proposition 5 of Merlev\`ede and Peligrad (2012) with $\varphi(u)=u^q$ and $x=\lambda  \psi(2^r) $, and using stationarity, we  derive that for any integer $r\ge 1$,
$$ {\mathbf P}\Big (\max_{1\leq i \le 2^r} |S_i-M_i|   \geq  \lambda  \psi(2^r) \Big )
 \ll   \frac{\|R_{2^r}\|_{q}^q}{\lambda^q ( \psi(2^r) )^{q}}
  +\frac{2^r}{\lambda^p ( \psi(2^r) )^{p} }\bigg(\sum_{l=0}^{r-1}2^{-l/p}\|\E (S_{2^l} |{ \mathcal F}_0)\|_{p}\bigg)^p \, .
$$
Notice now that by stationarity,  for all $i,j \geq 0$, $\| R_{i+j}\|_{q} \leq \| R_{i}\|_{q} + \| R_{j}\|_{q}$ and $\| \E (S_{i+j} |{ \mathcal F}_0)\|_{p} \leq \| \E (S_{i} |{ \mathcal F}_0)\|_{p}  + \| \E (S_{j} |{ \mathcal F}_0)\|_{p} $. Therefore applying Lemma \ref{claim1} of the appendix respectively with $V_n = ( \psi(n)  )^{-q}\| R_{n}\|^q_{q}$ and $V_n = \|\E (S_{n} |{ \mathcal F}_0)\|_{p}$, we infer that \eqref{but1thminter} will hold true
if \eqref{psintercond} is satisfied. $\diamond$

\subsection{Proof of Proposition \ref{approxrnq}}
Notice first that the following decomposition is valid:  for any positive integer $n$,
\beq \label{decRn*}
R_n  =  \sum_{k=1}^n \big ( X_k - \sum_{j =1}^n P_j(X_k) \big ) - \sum_{k=1}^n \sum_{j \geq n+1} P_k(X_j) = {\mathbf E} (S_n |{\mathcal F}_0) - \sum_{k=1}^n \sum_{j \geq n+1} P_k(X_j)
 \, .
\eeq
Let $N$ be a positive integer and write
\begin{eqnarray}  \label{decRn**}
\sum_{k=1}^n \sum_{j \geq n+1} P_k(X_j) & =& \sum_{j = n+1}^{n + N}  \sum_{k=1}^n P_k(X_j)  + \sum_{k=1}^n \sum_{j \geq n+N +1} P_k(X_j) \\
& = & \bkE (S_{n+N} - S_n |{\mathcal F}_{n}) - \bkE (S_{n+N} - S_n  |{\mathcal F}_{0}) +  \sum_{k=1}^n \sum_{j \geq n+N +1} P_k(X_j) \, .
\end{eqnarray}
Starting from \eqref{decRn*} and considering \eqref{decRn**}, item 1  follows. To prove item 2, we start from item 1 and use stationarity, to derive that for any positive integers $n$ and $N$,
\beq \label{normeprn*}
\Vert R_n \Vert_{p} \leq   \Vert{\mathbf E} (S_n |{\mathcal F}_0)  \Vert_p + 2 \Vert \bkE (S_{N}   |{\mathcal F}_{0}) \Vert_p +   \Big \Vert  \sum_{k=1}^n \sum_{j \geq n+N +1} P_k(X_j) \Big \Vert_{p} \, .
\eeq
Applying  then Burkholder's inequality and using the stationarity,
we obtain that for any positive integer $n$, there exists a positive constant $c_p$ such that
\begin{eqnarray} \label{Burkh1}
\Big \Vert  \sum_{k=1}^n \sum_{j \geq n+N +1} P_k(X_j) \Big \Vert^{p'}_{p}   \leq   c_p  \sum_{k=1}^n \Big \Vert \sum_{j \geq n+N +1} P_k(X_j) \Big \Vert_p^{p'}  = c_p  \sum_{k=1}^n  \Big \Vert  \sum_{j \geq N+k}  P_0(X_j) \Big \Vert_p^{p'} \, ,
\end{eqnarray}
where $p'=\min (2,p)$. Starting from (\ref{normeprn*}) and using (\ref{Burkh1}), item 2 follows. $\diamond$

\subsection{Proof of Theorem \ref{ThNAS}}

The proof will follow from Theorem \ref{interwunous} if we can show that under \eqref{Cond1cobSn},  \eqref{Cond1cobSn2} and  \eqref{condcarre}, the second part of condition \eqref{newcondTH22} is satisfied. With this aim, let $M_n=S_n-R_n$ and write that
$$
\Vert  \bkE (M_n^2 | {\mathcal F}_0  ) -\bkE  (M_n^2)
 \big \Vert_{p/2} \leq \Vert  \bkE (S_n^2 | {\mathcal F}_0  ) -\bkE  (S_n^2)
 \big \Vert_{p/2} +2 \Vert \bkE_0  (S_nR_n  ) -
\bkE  (S_nR_n)  \Vert _{p/2} + 2 \Vert R_n \Vert^{2}_{p} \, .
$$
Let $\beta_n= n^2 (\log n)^{(t-1)p/2}$. Since \eqref{condcarre} holds true, the second part of condition \eqref{newcondTH22} will be satisfied if
\beq \label{I*}
\sum_{n \geq 1} \frac{ 1 }{ \beta_n }\Vert R_n \Vert^{p}_{p} < \infty \, \text{ and } \sum_{n \geq 1} \frac{ 1 }{ \beta_n } \Vert \bkE_0  (S_nR_n  ) -
\bkE  (S_nR_n)  \Vert^{p/2}_{p/2} < \infty\, .\eeq
By using item 2 of Proposition \ref{approxrnq} with $N= u_n$, the first part of \eqref{I*} clearly holds under \eqref{Cond1cobSn}. To prove the second part of \eqref{I*}, we first notice that
$$ \Vert \bkE  (S_n \bkE  (S_n | {\mathcal F}_0  ) | {\mathcal F}_0  ) -
\bkE  (S_n {\mathbf E} (S_n |{\mathcal F}_0))  \Vert_{p/2}
  \leq 2 \Vert {\mathbf E} (S_n |{\mathcal F}_0)  \Vert^2_{p}   \, .
$$ Hence the first part of (\ref{Cond1cobSn}) implies that
\beq \label{p0II}
\sum_{n \geq 1} \frac{ 1 }{ \beta_n } \Vert \bkE^2  (S_n | {\mathcal F}_0  ) -
\bkE  (S_n {\mathbf E} (S_n |{\mathcal F}_0))  \Vert^{p/2}_{p/2} < \infty \, .
\eeq
In addition, $ \Vert \bkE  (S_n {\mathbf E} (S_{2n} - S_n |{\mathcal F}_0)| {\mathcal F}_0  ) -
\bkE  (S_n {\mathbf E} (S_{2n} - S_n |{\mathcal F}_0))  \Vert_{p/2} \leq  2 \Vert {\mathbf E} (S_n |{\mathcal F}_0)  \Vert^2_{p} $. Therefore, we also have that
\beq \label{p0bisII}
\sum_{n \geq 1}  \frac{1}{\beta_n}  \Vert \bkE  (S_n {\mathbf E} (S_{2n} - S_n |{\mathcal F}_0)| {\mathcal F}_0  ) -
\bkE  (S_n {\mathbf E} (S_{2n} - S_n |{\mathcal F}_0)) \Vert^{p/2}_{p/2} < \infty \, .
\eeq
Now since $S_n$ is ${\mathcal F}_n$-measurable, we get that
\begin{multline*} \Vert \bkE ( S_n {\mathbf E} (S_{2n}- S_n |{\mathcal F}_n)|{\mathcal F_0}) -\bkE ( S_n {\mathbf E} (S_{2n}- S_n |{\mathcal F}_n)) \Vert_{p/2} \\   = \Vert  {\mathbf E} (S_n(S_{2n}- S_n) |{\mathcal F_0}) - \bkE ( S_n (S_{2n}- S_n )) \Vert_{p/2} \, .
\end{multline*}
Next, using the identity $2ab = (a+b)^2 - a^2 -b^2$ and the stationarity, we obtain that
\begin{multline*}2 \Vert \bkE ( S_n {\mathbf E} (S_{2n}- S_n |{\mathcal F}_n)|{\mathcal F_0}) - \bkE ( S_n {\mathbf E} (S_{2n}- S_n |{\mathcal F}_n)) \Vert_{p/2} \\
 \quad \leq \Vert  {\mathbf E} (S^2_{2n}|{\mathcal F_0}) - \bkE ( S^2_{2n} ) \Vert_{p/2} + 2 \Vert  {\mathbf E} (S^2_{n}|{\mathcal F_0}) - \bkE ( S^2_{n} ) \Vert_{p/2} \, ,
\end{multline*}
which combined with (\ref{condcarre}) implies that
\beq \label{p1II}
\sum_{n \geq 1}  \frac{1}{\beta_n} \Vert \bkE ( S_n {\mathbf E} (S_{2n}- S_n |{\mathcal F}_n)|{\mathcal F_0}) - \bkE ( S_n {\mathbf E} (S_{2n}- S_n |{\mathcal F}_n)) \Vert^{p/2}_{p/2} < \infty \, .
\eeq
Therefore by combining (\ref{p0II}), (\ref{p0bisII}) and (\ref{p1II}), we derive that the  second part of (\ref{I*}) will be satisfied provided that
\beq \label{p2II}
\sum_{n \geq 1}  \frac{1}{\beta_n^{ 2/p} n^{1+2/p}} \Vert \bkE ( S_n {\widetilde R}_n |{\mathcal F_0})- \bkE (  S_n {\widetilde R}_n )\Vert_{p/2} < \infty \, ,
\eeq
where ${\widetilde R}_n = R_n -  {\mathbf E} (S_n |{\mathcal F}_0) - {\mathbf E} (S_{2n}- S_n |{\mathcal F}_n)  + {\mathbf E} (S_{2n} - S_n |{\mathcal F}_0)$. To prove (\ref{p2II}), we first write that
\begin{multline} \label{decp2II}
\Vert \bkE_0  (S_n{\widetilde R}_n  )  -
\bkE  (S_n{\widetilde R}_n )  \Vert_{p/2} \leq 2 \Vert \bkE_0  (S_n{\widetilde R}_n  ) \Vert_{p/2} \leq 2 \Vert \bkE^{1/2}_0  (S_n^2) \bkE^{1/2}_0({\widetilde R}_n ^2  ) \Vert_{p/2} \\
 \leq  2 \Vert ( \bkE_0  (S_n^2)- \bkE(S_n^2) )^{1/2}\bkE^{1/2}_0({\widetilde R}^2_n  ) \Vert_{p/2}+ 2 (\bkE(S_n^2) )^{1/2}\Vert \bkE^{1/2}_0({\widetilde R}_n^2 ) \Vert_{p/2}  \\
 \leq   \Vert \bkE_0  (S_n^2)- \bkE(S_n^2)  \Vert_{p/2} + \Vert {\widetilde R}_n \Vert^2_{p} + 2 (\bkE(S_n^2) )^{1/2}\Vert \bkE^{1/2}_0({\widetilde R}_n^2  ) \Vert_{p/2} \, .
\end{multline}
In addition,
\beq \label{rntilde}
\Vert {\widetilde R}_n \Vert_{p} \leq  \Vert R_n \Vert_{p} + 3 \Vert  {\mathbf E} (S_n |{\mathcal F}_0) \Vert_p  \ \text{ and } \ (\bkE(S_n^2) )^{1/2}\Vert \bkE^{1/2}_0({\widetilde R}_n^2  ) \Vert_{p/2}  \ll  n^{1/2} \Vert {\widetilde R}_n  \Vert_2\, ,
\eeq
where for the last inequality we have used the fact that the function $x \mapsto |x|^{p/4}$ is concave and that $\bkE(S_n^2) \ll n$. Next, taking into account
item 1 of Proposition \ref{approxrnq}   with $N=n$, the following inequality holds:
\beq \label{rntilde2}
\Vert {\widetilde R}_n  \Vert^2_2 =  \sum_{k=1}^n \Big\| \sum_{j \geq k+n}P_0 (X_j) \Big\|^2_2 \, ,
\eeq
Therefore, starting from (\ref{decp2II}) and using (\ref{rntilde}) and (\ref{rntilde2}), we then infer that (\ref{p2II}) holds by taking into account (\ref{Cond1cobSn2}), (\ref{condcarre}) and the first part of (\ref{I*}). This ends the proof of  the second part of (\ref{I*}) and then of the theorem. $\diamond$

\subsection{Proof of Proposition \ref{directcond}}

We first notice that  the first part of (\ref{Cond1cob*}) implies that 
$\|{\mathbf E}(X_n | {\mathcal F}_0) \|_p=o(n^{2/p^2-1} (\log n)^{(t-1)/2})$, so that 
$
\sum_{n>0} n^{-1/p} \|{\mathbf E}(X_n | {\mathcal F}_0) \|_p < \infty \, ,
$
since $p>2$.  
This  implies that 
$X_0$ is regular in the sense that ${\mathbf E}(X_0 |{\mathcal F}_{- \infty})=0$ almost surely.
Applying  Lemma \ref{complemmap0xi} with $q=1$, it follows that 
$ 
\sum_{n>0} \|P_0(X_n)\|_p < \infty
$
and the condition (\ref{Gor}) is satisfied.
Now, since $X_0$ is regular, 
the condition  $\sum_{n>0} \|P_0(X_n)\|_2 < \infty$ implies that
the series $ \sigma^2 =\sum_{k \in {\mathbf Z}} {\rm Cov } (X_0 , X_k)$ converges absolutely, 
and consequently $n^{-1} {\mathbf E}(S_n^2)$ converges to  $\sigma^2$.

Let $\beta_n= n^2 (\log n)^{(t-1)p/2}$. Let us prove that  \eqref{Cond1cobSn} holds with $u_n=[n^{p/2}]$ as soon as the first part of \eqref{Cond1cob*} is satisfied. Since $\max_{k \in \{n, u_n \}}\| {\mathbf E}(S_{k} | {\mathcal F}_0) \|_p \leq \sum_{k=1}^{u_n}\| {\mathbf E}(X_k| {\mathcal F}_0) \|_p$, H\"older's inequality implies that for any $\gamma \in ]0, 2/p[$,
\begin{align*}
\sum_{n \geq 2}  \frac{ 1 }{ \beta_n}\max_{k \in \{n, u_n \}}\| {\mathbf E}(S_{k} | {\mathcal F}_0) \|_p^p&  \ll \sum_{n \geq 2}  \frac{ n^{ \gamma p/2} }{  \beta_n} \sum_{k=1}^{[n^{p/2}]} k^{p-1- \gamma}\| {\mathbf E}(X_k| {\mathcal F}_0) \|^p_p \, .
\end{align*}
Therefore changing the order of summation, we infer that  the first part of \eqref{Cond1cobSn} holds with  $u_n=[n^{p/2}]$ as soon as the first part of \eqref{Cond1cob*} does. Let us prove now that the second part of \eqref{Cond1cobSn} is satisfied with $u_n=[n^{p/2}]$. Using Lemma \ref{complemmap0xi} given in Appendix, we first get that
$$
\sum_{n \geq 2}   \frac{ 1 }{ \beta_n }  \Big ( \sum_{k=1}^n \Big \| \sum_{j \geq k+u_n}P_0 (X_j) \Big \|^2_p \Big )^{p/2} \ll \sum_{n \geq 2}   \frac{ n^{p/2} }{ \beta_n }  \Bigr ( \sum_{k \geq [n^{p/2}/2]  } \frac{\Vert \E(X_k |{\mathcal F}_0 )\Vert_p}{k^{1/p}} \Bigl )^{p} \, .
$$
Next H\"older's inequality implies that for any $\gamma \in ]1-1/p, 1-2/p^2[$,
$$
\sum_{n \geq 2}   \frac{ 1 }{ \beta_n }  \Big ( \sum_{k=1}^n \Big \| \sum_{j \geq k+u_n}P_0 (X_j) \Big \|^2_p \Big )^{p/2} \ll \sum_{n \geq 2}   \frac{ n^{p( p - \gamma p)/2} }{ \beta_n }   \sum_{k \geq [n^{p/2}]  } \frac{k^{\gamma p}\Vert \E(X_k |{\mathcal F}_0 )\Vert^p_p}{k}  \, .
$$
Changing the order of summation, we then derive that the second part of \eqref{Cond1cobSn} holds with  $u_n=[n^{p/2}]$ as soon as the first part of \eqref{Cond1cob*} does.

It remains to prove that the second part of \eqref{Cond1cob*} entails that \eqref{Cond1cobSn2} holds. Using again Lemma \ref{complemmap0xi} followed by an application of H\"older's inequality, we first obtain that for any $\gamma \in ]1-2/p,2-4/p  [$,
\begin{align*}
\sum_{n \geq 2}   \frac{ n^{p/4} }{ \beta_n }
  \Big ( \sum_{k=1}^n \Big \| \sum_{j \geq k+n}P_0 (X_j) \Big \|^2_2 \Big )^{p/4} & \ll \sum_{n \geq 2}   \frac{ n^{p/2} }{ \beta_n } \Bigr ( \sum_{k \geq [n/2]  } \frac{\Vert \E(X_k |{\mathcal F}_0 )\Vert_2}{k^{1/2}} \Bigl )^{p/2} \\
  & \ll \sum_{n \geq 2}   \frac{ n^{p( 2 -  \gamma -2/p)/2 }  }{ \beta_n }  \sum_{k \geq n}  k^{\gamma p/2} \frac{\Vert \E(X_k |{\mathcal F}_0 )\Vert^{p/2}_2}{k^{p/4}} \, .
\end{align*}
Changing the order of summation, it follows that  \eqref{Cond1cobSn2} holds  provided that the second part of \eqref{Cond1cob*} does. This ends the proof of the proposition. $\diamond$

\subsection{Proof of Corollary \ref{coralphaphi}.} Notice first that, since $\gamma \in ]0,1]$ and $p >2$, the condition (\ref{cond1coralphaphi}) implies the second part of (\ref{Cond1cob*}). Now, since $( \Vert \bkE  (X_n |{\mathcal F}_0 )\Vert_p)_{n \geq 1}$ is a  nonincreasing sequence,  condition  (\ref{cond1coralphaphi}) implies that \beq \label{easy} \Vert \bkE  (X_n |{\mathcal F}_0 )\Vert^p_p = o \big ((\log n)^{(t-1) p} n^{3 -2p} \big ) \, ,
\eeq which in turn entails that  the first part of (\ref{Cond1cob*})
is satisfied. By Proposition \ref{directcond}, it follows that (\ref{Gor}), (\ref{Cond1cobSn}) and (\ref{Cond1cobSn2}) are satisfied.

It remains to prove that condition (\ref{condcarre}) holds as soon as (\ref{cond1coralphaphi}) and (\ref{cond2coralphaphi}) hold.  With this aim, setting
$$
\gamma({\mathcal F}_0, m,k) := \Vert \bkE (X_mX_{k+m} |
{\mathcal F}_0) - \bkE (X_mX_{k+m}) \Vert_{p/2} \, ,
$$
we first write that for any $\gamma \in [0,1]$,
\begin{align} \label{corphidec1}
 \Vert \bkE (S_n^2 | {\mathcal F}_0)  - \bkE (S_n^2) \Vert_{p/2}
 & \leq 2 \sum_{m=1}^n\sum_{k=0}^{n-m} \gamma({\mathcal F}_0, m,k) \nonumber \\
& \leq 2 \sum_{m=1}^n\sum_{k=0}^{[m^{\gamma}]} \gamma({\mathcal F}_0, m,k) + 2 \sum_{m=1}^n\sum_{k=[m^{\gamma}]+1}^n \gamma({\mathcal F}_0, m,k)\, .
\end{align}
We shall bound up $\gamma({\mathcal F}_0, m,k) $ in two ways. First  we consider the bound: for any $k\geq 0$,
$$
\gamma({\mathcal F}_0, m,k)  \leq \sup_{i \geq j \geq m }\| {\mathbf E}(X_iX_j | {\mathcal F}_0) - {\mathbf E}(X_iX_j)\|_{p/2} = \widetilde \gamma (m)\, .
$$
Therefore, by H\"older's inequality,
\begin{align*}
\sum_{n \geq 2} \frac{(\log n)^{(1-t)p/2}}{n^2}&  \Big ( \sum_{m=1}^n\sum_{k=0}^{[m^{\gamma}]} \gamma({\mathcal F}_0, m,k)\Big )^{p/2}
 \ll \sum_{n \geq 2} \frac{(\log n)^{(1-t)p/2}}{n^{2}} \Big ( \sum_{m=1}^n m^{\gamma}\widetilde \gamma (m)\Big )^{p/2} \\
& \ll \sum_{n \geq 2} \frac{n^{(1 - \alpha)(p-2)/2}(\log n)^{(1-t)p/2}}{n^{2}} \sum_{m=1}^n m^{\gamma p/2 + \alpha(p-2)/2 }(\widetilde \gamma (m) )^{p/2} \, ,
\end{align*}
for any $\alpha \in ] 1 -2/(p-2), 1[$. Changing the order of summation, it follows that under \eqref{cond2coralphaphi},
\beq \label{b11rst}
\sum_{n \geq 2} \frac{(\log n)^{(1-t)p/2}}{n^{2}}   \Big ( \sum_{m=1}^n\sum_{k=0}^{[m^{\gamma}]} \gamma({\mathcal F}_0, m,k)\Big )^{p/2}
 \leq \sum_{n \geq 2} \frac{(\log n)^{(1-t)p/2}n^{(1+\gamma)p/2}}{n^{2}} (\widetilde \gamma (n) )^{p/2}< \infty\, .
\eeq
We write now that
\begin{multline} \label{b12nd}
\gamma ({\mathcal F}_0, m,k)  \leq 2 \|{\mathbf E}(X_mX_{k+m}|{\mathcal F}_0)\|_{p/2}
\\ \leq  2  \|{\mathbf E}\big ( (X_m - \bkE_0(X_m) ) (X_{k+m} - \bkE_0 ( X_{k+m}))|{\mathcal F}_0 \big )\|_{p/2} + 2 \|\bkE_0(X_m)\bkE_0 ( X_{k+m} )\|_{p/2} \, .
\end{multline}
We notice first that H\"older's inequality entails that
\begin{align*}
\sum_{n \geq 2} \frac{(\log n)^{(1-t)p/2}}{n^{2}}  & \Big ( \sum_{m=1}^n\sum_{k=0}^{n-m} \|\bkE_0(X_m)\bkE_0 ( X_{k+m} )\|_{p/2} \Big )^{p/2}
  \ll \sum_{n \geq 2} \frac{(\log n)^{(1-t)p/2}}{n^{2}} \Big ( \sum_{k=1}^{n} \|\bkE_0(X_k)\|_p \Big )^{p}
\\
& \ll \sum_{n \geq 2} \frac{(\log n)^{(1-t)p/2}  \, n^{(1- \alpha) (p-1)}}{n^{2}} \sum_{k=1}^{n}k^{ \alpha (p-1)} \|\bkE_0(X_k)\|_p^{p} \, ,
\end{align*}
for any $\alpha \in ] 1 -1/(p-1), 1[$. Changing the order of summation and using (\ref{easy}), it follows that
\begin{align} \label{b12nd*}
\sum_{n \geq 2} \frac{(\log n)^{(1-t)p/2}}{n^{2}}    \Big ( \sum_{m=1}^n\sum_{k=0}^{n-m} & \|\bkE_0(X_m)\bkE_0 ( X_{k+m} )\|_{p/2} \Big )^{p/2}  \nonumber \\ & \ll \sum_{n \geq 2} \frac{(\log n)^{(1-t)p/2}  \, n^{p}}{n^{2}}  \|\bkE_0(X_n)\|_p^{p}  < \infty \, .
\end{align}
Starting from \eqref{corphidec1} and using \eqref{b11rst}, \eqref{b12nd} and \eqref{b12nd*}, condition (\ref{condcarre}) will be satisfied if we can prove that
\beq \label{b1third}
\sum_{n \geq 2} \frac{(\log n)^{(1-t)p/2}}{n^{2}}    \Big ( \sum_{m=1}^n\sum_{k=[m^{\gamma}]+1}^n \gamma^*(m,k) \Big )^{p/2}    < \infty \, .
\eeq
where $\gamma^*(m,k):= \|{\mathbf E} \big ( (X_m - \bkE_0(X_m) ) (X_{k+m} - \bkE_0 ( X_{k+m}))|{\mathcal F}_0 \big )\|_{p/2} $. With this aim, since $ X_m - \bkE_0(X_m)  = \sum_{\ell=1}^m P_{\ell}( X_m)$, we first observe that
$$
 \gamma^*(m,k) \leq
 \Big \Vert {\mathbf E}\big (\sum_{\ell=1}^m
 P_{\ell} (X_m)  P_{\ell}( X_{k+m})|{\mathcal F}_0 \big )\Big \Vert_{p/2}  \leq  \sum_{\ell=1}^m \Vert P_{\ell} (X_m)\Vert_p \Vert P_{\ell} (X_{k+m})\Vert_p\, .
$$
Therefore by stationarity,
\beq \label{b12nd*bis}
\gamma^*(m,k) \leq  \sum_{\ell=0}^{m-1} \Vert P_{0} (X_{\ell})\Vert_p \Vert P_{0} (X_{\ell +k})\Vert_p\, .
\eeq
Now we write that
\beq \label{b12nd*ter}
\sum_{m=1}^n\sum_{k=[m^{\gamma}] +1}^n \gamma^*(m,k)  \ll
 \sum_{k=1}^{[n^{\gamma}]}\sum_{m=1}^{[k^{1/\gamma}] +1 } \gamma^*(m,k) +  \sum_{k=[n^{\gamma}]+1}^n\sum_{m=1}^{n} \gamma^*(m,k) \, .
\eeq
Using \eqref{b12nd*bis} and H\"older's inequality twice, we first derive that
\begin{align*}
\Big ( \sum_{k=1}^{[n^{\gamma}]}& \sum_{m=1}^{[k^{1/\gamma}] +1 } \gamma^*(m,k)\Big )^{p/2}  \ll \Big (\sum_{\ell=0}^{n-1} \Vert P_{0} (X_{\ell})\Vert_p \sum_{k=1}^{[n^{\gamma}]}k^{1/\gamma} \ \Vert P_{0} (X_{\ell +k})\Vert_p\Big )^{p/2} \\
& \ll \Big (\sum_{i=0}^{n-1} \Vert P_{0} (X_{i})\Vert_p\Big )^{\frac{p-2}{2}} \  \sum_{\ell=0}^{n-1} \Vert P_{0} (X_{\ell})\Vert_p  \Big ( \sum_{k=1}^{[n^{\gamma}]}k^{1/\gamma} \ \Vert P_{0} (X_{\ell +k})\Vert_p \Big )^{p/2} \\
& \ll n^{\gamma(1- \alpha) (p-2)/2}\Big (\sum_{i=0}^{n-1} \Vert P_{0} (X_{i})\Vert_p\Big )^{\frac{p-2}{2}} \ \sum_{\ell=0}^{n-1} \Vert P_{0} (X_{\ell})\Vert_p   \sum_{k=1}^{[n^{\gamma}]}k^{p/(2\gamma)} k^{\alpha(p-2)/2}\ \Vert P_{0} (X_{\ell +k})\Vert_p^{p/2}\, .
\end{align*}
for any $\alpha \in ] 1 -2/(p-2), 1[$. We then infer that
\begin{align*}
\sum_{n \geq 2}  \frac{(\log n)^{(1-t)p/2}}{n^{2}} & \Big ( \sum_{k=1}^{[n^{\gamma}]} \sum_{m=1}^{[k^{1/\gamma}] +1 } \gamma^*(m,k)\Big )^{p/2}  \\
& \ll \Big (\sum_{\ell\geq 0} \Vert P_{0} (X_{\ell})\Vert_p\Big )^{p/2}  \sum_{n \geq 2}  \frac{(\log n)^{(1-t)p/2} \, n^{p/(2\gamma) +p/2}}{n^{1+ 1/\gamma}}\ \Vert P_{0} (X_{n})\Vert_p^{p/2}\, .
\end{align*}
Notice now that by Lemma \ref{complemmap0xi}, $\sum_{\ell\geq 0} \Vert P_{0} (X_{\ell})\Vert_p < \infty$ as soon as $\sum_{k >0} k^{-1/p}\Vert \E_{0} (X_{\ell})\Vert_p < \infty$, which clearly holds under \eqref{easy} since $p>2$.  On an other hand, using Lemma \ref{complemmap0xi} again, we get that
\begin{align} \label{bfourth}
\sum_{n \geq 2}  & \frac{(\log n)^{(1-t)p/2} \, n^{p/(2\gamma) +p/2}} {n^{1 + 1/\gamma}} \ \Vert P_{0} (X_{n})\Vert_p^{p/2}\ll \sum_{n \geq 2}   \frac{(\log n)^{(1-t)p/2} \, n^{p/(2\gamma) +p/2}} {n^{2 + 1/\gamma}} \sum_{k \geq n} \Vert P_{0} (X_{k})\Vert_p^{p/2} \nonumber \\
& \ll \sum_{n \geq 2}   \frac{(\log n)^{(1-t)p/2} \, n^{p/(2\gamma) +p/2}} {n^{2 + 1/\gamma}} \sum_{k \geq [n/2]} k^{-1/2}\Vert \E_{0} (X_{k})\Vert_p^{p/2}\, .
\end{align}
It follows that under (\ref{cond1coralphaphi}),
$$
\sum_{n \geq 2}  \frac{(\log n)^{(1-t)p/2}}{n^{2}}  \Big ( \sum_{k=1}^{[n^{\gamma}]} \sum_{m=1}^{[k^{1/\gamma}] +1 } \gamma^*(m,k)\Big )^{p/2} < \infty \, .
$$
Therefore from \eqref{b12nd*ter}, \eqref{b1third} will follow if we can prove that
\beq \label{b15thbis}
\sum_{n \geq 2}  \frac{(\log n)^{(1-t)p/2}}{n^{2}} \Big (  \sum_{k=[n^{\gamma}]+1}^n  \sum_{m=1}^{n} \gamma^*(m,k)\Big )^{p/2}< \infty \, .
\eeq
Using \eqref{b12nd*bis} and H\"older's inequality twice,
\begin{align*}
& \Big (  \sum_{k=[n^{\gamma}]+1}^n  \sum_{m=1}^{n} \gamma^*(m,k)\Big )^{p/2}  \ll n^{p/2} \Big (\sum_{\ell=0}^{n-1} \Vert P_{0} (X_{\ell})\Vert_p \sum_{k=[n^{\gamma}]+1}^n   \Vert P_{0} (X_{\ell +k})\Vert_p\Big )^{p/2} \\
& \ll n^{p/2}  \Big (\sum_{i=0}^{n-1} \Vert P_{0} (X_{i})\Vert_p\Big )^{(p-2)/2} \sum_{\ell=0}^{n-1} \Vert P_{0} (X_{\ell})\Vert_p  \Big ( \sum_{k=[n^{\gamma}]+1}^n \Vert P_{0} (X_{\ell +k})\Vert_p \Big )^{p/2} \\
& \ll n^{p/2} n^{\gamma(1- \alpha) (p-2)/2}\Big (\sum_{i=0}^{n-1} \Vert P_{0} (X_{i})\Vert_p\Big )^{(p-2)/2} \sum_{\ell=0}^{n-1} \Vert P_{0} (X_{\ell})\Vert_p    \sum_{k=[n^{\gamma}]+1}^n k^{\alpha(p-2)/2}\ \Vert P_{0} (X_{\ell +k})\Vert_p^{p/2}\, ,
\end{align*}
for any $\alpha \in ] 1, 1 + 1/ \gamma[$. It follows that
\begin{align*}
\sum_{n \geq 2}  \frac{(\log n)^{(1-t)p/2}}{n^{2}} &\Big (  \sum_{k=[n^{\gamma}]+1}^n  \sum_{m=1}^{n} \gamma^*(m,k)\Big )^{p/2} \\
& \ll \Big (\sum_{\ell\geq 0} \Vert P_{0} (X_{\ell})\Vert_p\Big )^{p/2}  \sum_{n \geq 2}  \frac{(\log n)^{(1-t)p/2} \,
n^{p/(2\gamma) +p/2}}{n^{1 + 1/\gamma}}\ \Vert P_{0} (X_{n})\Vert_p^{p/2}\, ,
\end{align*}
which is finite under (\ref{cond1coralphaphi}) (see \eqref{bfourth}). This ends the proof of Corollary \ref{coralphaphi}. $\diamond$

\subsection{Proof of Corollary \ref{coralpha}.}
Starting from \eqref{corphidec1} with $\gamma=1$ combined with the bound
$$
\gamma({\mathcal F}_0, m,k)  \leq  \Big ( \sup_{i \geq j \geq m }\| {\mathbf E}(X_iX_j | {\mathcal F}_0) - {\mathbf E}(X_iX_j)\|_{p/2} \Big )\wedge 
\Big (2 \|X_0 {\mathbf E}(X_k | {\mathcal F}_0) \|_{p/2} \Big ) \, ,
$$
where $\gamma({\mathcal F}_0, m,k) $ is defined in (\ref{corphidec1}), we infer that (\ref{condcarre}) is satisfied provided that (\ref{cond2coralpha}) and (\ref{cond2coralpha*}) are. The corollary then follows from an application of Theorem \ref{ThNAS} together with Proposition \ref{directcond}. $\diamond$

\subsection{Proof of Corollary \ref{corMarkovChain}.}
Let us first check that (\ref{Gor}) is satisfied. 
By  definition of the  transition probability, 
and since $f$ is an odd function and $\upsilon$ is symmetric,
\begin{equation}
{\mathbf{E}}_{0}(X_{n})= {\mathbf{E}}(f(\zeta_{n})|\zeta_{0})=(1-|\zeta_{0}|)^{n}f(\zeta_{0})\text{
a.s.} \label{operator}%
\end{equation}
Taking into account that $|f(x)|\leq C|x|^{1/2}$ and  the assumption on $\upsilon$
it follows that 
$$
  \|{\mathbf{E}}_{0}(X_{n})\|_p^p \ll \int_0^1 (1-x)^{np} x^{p/2} x^{a-1} dx \, ,
$$
and  using the properties of the Beta functions,
$
\|{\mathbf{E}}_{0}(X_{n})\|_p^p =O(n^{-a-p/2})
$. Hence,  if $a>(p-2)/2$,
$\sum_{n>0} n^{-1/p} \|{\mathbf{E}}_{0}(X_{n})\|_p < \infty$ and by Lemma 
\ref{complemmap0xi}
of the appendix, $\sum_{n>0} \|P_0(X_n)\|_p < \infty$, proving that (\ref{Gor}) is 
satisfied. We also get that 
\beq \label{rio1}
\Vert {\mathbf{E}}_{0}(S_{n}(f))\Vert_{p}\ll \max( \log n , n^{1/2-a/p} ) \, .
\eeq

We prove now that for $M_n=\sum_{k=1}^n d_i$ where $d_i= \sum_{n \geq i}P_i(X_n)$,
\beq \label{estimatemn}
\Vert {\mathbf{E}}_{0}(M^2_{n})-{\mathbf{E}}(M^2_{n}) \Vert_{p/2} \ll \max( \log n , n^{1-2a/p} )  \, .
\eeq
Starting from (\ref{operator}), we infer that
$$
d_1= \frac{f(\zeta_{1}) }{|\zeta_{1}|} - \frac{f(\zeta_{0}) }{|\zeta_{0}|} + f(\zeta_{0}) \, \text{
a.s.}  \, ,
$$
yielding to
$$
M_n = \frac{f(\zeta_{n}) }{|\zeta_{n}|} ( 1 - |\zeta_n |)  - \frac{f(\zeta_{0}) }{|\zeta_{0}|} ( 1 - |\zeta_0 |)  + S_n(f) \, \text{
a.s.}  \, .
$$
Hence \eqref{estimatemn} will be proven if we can show that
\beq \label{irre1}
\Vert {\mathbf{E}}_{0}(S^2_{n}(f))-{\mathbf{E}}(S^2_{n}(f)) \Vert_{p/2} \ll \max( \log n , n^{1-2a/p} ) \, ,
\eeq
\beq \label{irre2}
\Big \Vert {\mathbf{E}}_{0}\Big ( \big ( f(\zeta_{n}) (  |\zeta_{n}|^{-1} - 1 ) -  f(\zeta_{0}) ( |\zeta_{0}|^{-1} - 1 )  \big )^2 \Big )  \Big \Vert_{p/2} \ll \max( \log n , n^{1-2a/p} )\, ,
\eeq
and
\beq \label{irre3}
\Big \Vert {\mathbf{E}}_{0}\Big ( S_n (f) \big ( f(\zeta_{n}) (  |\zeta_{n}|^{-1} - 1 ) -  f(\zeta_{0}) ( |\zeta_{0}|^{-1} - 1 )  \big ) \Big )  \Big \Vert_{p/2}  \ll \max( \log n , n^{1-2a/p} ) \, .
\eeq
The bound \eqref{irre1} has been proved in Rio (2009) (see the computations leading to his bound (4.15)). We turn now to the proof of (\ref{irre2}). 
According to the definition of the transition probability, we first notice that for any positive integer $n$, and any function $g$ on $E$,
\beq \label{irre4}
{\mathbf{E}}_{0}\big ( g(\zeta_{n})    \big )= {\mathbf{E}}_{0}\Big ( (  1 - |\zeta_{n-1}| )g(\zeta_{n-1}) \Big )    + {\mathbf{E}}_{0} ( | \zeta_{n-1}|) \int_E g(y)  \upsilon (dy) \text{
a.s} \, ,
\eeq
and since $f$ is an odd function and $\upsilon$ is symmetric,
\beq \label{irre5}
{\mathbf{E}}_{0}\big ( f(\zeta_{n})  (  |\zeta_{n}|^{-1} - 1 ) \big )
  =    (1-|\zeta_{0}|)^{n+1} \frac{f(\zeta_{0})}{|\zeta_{0}|}  \text{
a.s} \, .
\eeq
Therefore, using \eqref{irre4} with $g(x) = y^{-2} (1-|y|)^2 f^2(y)$ and \eqref{irre5}, we get that for any positive integer $n$,
\begin{multline*}
 \Big  \Vert {\mathbf{E}}_{0}\Big ( \big ( f(\zeta_{n}) (  |\zeta_{n}|^{-1} - 1 ) -  f(\zeta_{0}) ( |\zeta_{0}|^{-1} - 1 )  \big )^2 \Big )  \Big \Vert_{p/2} \ll    \Vert  \zeta_{n-1} \Vert_{p/2} \int_E  (1-|y|)^2 \frac{f^2(y) }{y^2}d \upsilon (y) \\
 + \Big \Vert  {\mathbf{E}}_{0}\Big ( \frac{f^2(\zeta_{n-1}) }{\zeta_{n-1}^2}(  1 - |\zeta_{n-1}| )^3 -2 \frac{f^2(\zeta_{0}) }{\zeta_{0}^2}(  1 - |\zeta_{0}| )^{n+2}  + \frac{f^2(\zeta_{0}) }{\zeta_{0}^2}(  1 - |\zeta_{0}| )^2  \Big )   \Big \Vert_{p/2} \, .
\end{multline*}
With the help of \eqref{irre4}, we then derive after $n-1$ steps that
\begin{multline} \label{irre6}
 \Big  \Vert {\mathbf{E}}_{0}\Big ( \big ( f(\zeta_{n}) (  |\zeta_{n}|^{-1} - 1 ) -  f(\zeta_{0}) ( |\zeta_{0}|^{-1} - 1 )  \big )^2 \Big )  \Big \Vert_{p/2} \ll   \Vert  \zeta_{0} \Vert_{p/2}  \sum_{k=0}^{n-1}\int_E  (1-|y|)^{2+k} \frac{f^2(y) }{y^2}d \upsilon (y)  \\
+ \Big \Vert   \frac{f^2(\zeta_{0}) }{\zeta_{0}^2}(  1 - |\zeta_{0}| )^2 \big ( 1  -(  1 - |\zeta_{0}| )^{n} \big )     \Big \Vert_{p/2} \, .
\end{multline}
Taking into account that $|f(x)|\leq C|x|^{1/2}$, the assumption on $\upsilon$ and using the properties of the Beta functions, we get that for any positive integer $n$
\beq \label{irre7} \sum_{k=0}^{n-1}\int_E  (1-|y|)^{2+k} \frac{f^2(y) }{y^2}d \upsilon (y) \ll \max( \log n , n^{1-a} )  \, .
\eeq
Since $|f(x)|\leq C|x|^{1/2}$, by taking into account \eqref{hypopi}, it follows that for any positive integer $n$,
$$
 \Big  \Vert    \frac{f^2(\zeta_{0}) }{\zeta_{0}^2}(  1 - |\zeta_{0}| )^2 \big (1- (  1 - |\zeta_{0}| )^n   \big )   \Big \Vert^{p/2}_{p/2} 
 \ll  \int_0^1 x^{-p/2} (  1 - x )^{p} \big (  1 - (  1 - x )^{n} \big )^{p/2}   x^{a-1} dx  \, .
$$
Using the fact that $1 - (  1 - x )^{n} \leq \min ( 1, n x)$, 
we then infer that for any positive integer $n$,
\beq \label{irre8} \Big  \Vert   \frac{f^2(\zeta_{0}) }{\zeta_{0}^2}(  1 - |\zeta_{0}| )^2 \big ( 1-(  1 - |\zeta_{0}| )^n \big )      \Big \Vert^{p/2}_{p/2}  \ll \max( \log n, n^{p/2 -a }) \, .
\eeq
The inequality \eqref{irre6} combined with the bounds \eqref{irre7} and \eqref{irre8} gives \eqref{irre2}. It remains to show that (\ref{irre3}) holds true. Taking into account the definition of the transition probability, the fact that     $f$ is an odd function and $\upsilon$ is symmetric, it follows that
\beq \label{irre9}
{\mathbf{E}}_{0}\big ( S_n (f)f(\zeta_{n}) (  |\zeta_{n}|^{-1} - 1 )  \big ) = \sum_{k=1}^{n} {\mathbf{E}}_{0}\Big (  \frac{f^2(\zeta_{k}) }{|\zeta_{k}| }(  1 -|\zeta_{k}| )^{n-k+1}  \Big ) \, .
\eeq
Next, by using \eqref{irre4} $k$ times, we infer that for any positive integer $k$
\begin{multline} \label{irre10}
\Big  \Vert  {\mathbf{E}}_{0}\Big (   \frac{f^2(\zeta_{k}) }{|\zeta_{k}| }(  1 -|\zeta_{k}| )^{n-k+1}  \Big )\Big \Vert_{p/2} \leq \Vert  \zeta_{0} \Vert_{p/2}  \sum_{\ell=0}^{k-1}\int_E  (1-|y|)^{n-\ell} \frac{f^2(y) }{|y|}d \upsilon (y) \\
+ \big  \Vert f^2(\zeta_{0}) |\zeta_{0}|^{-1}(  1 -|\zeta_{0}| )^{n+1}  \big \Vert_{p/2} \, .
\end{multline}
On the other hand, since $f$ is an odd function and $\upsilon$ is symmetric, by the definition of the transition probability, it follows that
\beq \label{irre11}{\mathbf{E}}_{0}\big ( S_n (f)  f(\zeta_{0}) ( |\zeta_{0}|^{-1} - 1 ) \big ) = \frac{f^2(\zeta_{0}) }{|\zeta_{0}| } \sum_{k=1}^{n}  (  1 -|\zeta_{0}| )^{k+1}   \, .
\eeq
Therefore, taking into account \eqref{irre9}, \eqref{irre10} and \eqref{irre11} together with \eqref{hypopi} and the fact $|f(x)|\leq C|x|^{1/2}$, we infer that
\begin{multline*} \label{irre12}
\Big \Vert {\mathbf{E}}_{0}  \Big ( S_n (f) \big ( f(\zeta_{n}) (  |\zeta_{n}|^{-1} - 1 ) -  f(\zeta_{0}) ( |\zeta_{0}|^{-1} - 1 )  \big ) \Big )  \Big \Vert_{p/2} \ll \sum_{k=1}^n \sum_{\ell =k}^n\int_0^1  (1-x)^{\ell} x^a dx \nonumber\\
 + n  \Big ( \int_0^1  (1-x)^{np/2} x^{a-1} dx  \Big )^{2/p} + \sum_{k=1}^n \Big ( \int_0^1  (1-x)^{kp/2} x^{a-1} dx  \Big )^{2/p} \, .
\end{multline*}
The bound (\ref{irre3}) then follows from the properties of the Beta functions. This ends to proof of \eqref{estimatemn}.

It remains to use Theorem \ref{interwunous} combined with the bounds \eqref{rio1} and \eqref{estimatemn} and the assumption that $a>(p-2)/2$, to end the proof of the corollary. Indeed when $p=4$, this directly implies that $S_n(f)$ satisfies the strong approximation \eqref{strassen} with $\sigma^2 = \sum_{k \in {\mathbf Z}} {\rm Cov } (X_0 , X_k)$ and rate $b_n=n^{1/4} (\log n)^{3/4 + \varepsilon}$ for any $\varepsilon >0$. Now when $p \in ]2,4[$, since $a>(p-2)/2$, the arguments also hold for some $p' >p$ leading to the rate $b_n = n^{1/p}$ in \eqref{strassen}.  $\diamond$

\subsection{Proof of Corollaries \ref{circle0} and \ref{circle}.}
To prove Corollary \ref{circle0}, it suffices to prove
that the sequence $X_i= f(\xi_i)-m(f)$ satisfies the conditions
(\ref{cond1coralphaphi}) and (\ref{cond2coralphaphi}) of Corollary \ref{coralphaphi} with $t=1$.

Note that, according to the proofs of Lemmas 5.2 and 5.3 in Dedecker and Rio (2008) and to their inequality (5.18),
$$\|{\mathbf{E}}(X_n | \xi_0)\|_\infty \leq C_1 (f)\sum_{k \in {\mathbf{Z}}^*} |k|^{-s} \big ( \log ( 1 + |k|) \big )^{-(1+\varepsilon)}
|\cos(2\pi k a)|^n \, , $$
and
$$
\sup_{i \geq j \geq n }\| {\mathbf E}(X_iX_j | {\mathcal F}_0) - {\mathbf E}(X_iX_j)\|_{\infty} \leq C_2 (f)\sum_{k \in {\mathbf{Z}}^*} |k|^{-s} \big ( \log ( 1 + |k|) \big )^{-(1+\varepsilon)}
|\cos(2\pi k a)|^n  \, .
$$
Hence to verify the conditions (\ref{cond1coralphaphi}) and (\ref{cond2coralphaphi}), it suffices to take $\gamma$ such that $2/p - \gamma= 1/p - \gamma^{-1} ( 1- 2/p)$ that is
\beq \label{egalitegamme}
\gamma= \frac{1+ \sqrt{1+4p(p-2)}}{2p} \, ,
\eeq
and to  show that
\begin{equation} \label{aimcircle}
\sum_{n \geq 1} n^{\gamma - 2/p }\sum_{k \in {\mathbf{Z}}^*} |k|^{-s} \big ( \log ( 1 + |k|) \big )^{-(1+\varepsilon)}
|\cos(2\pi k a)|^n  < \infty \, .
\end{equation} Note first that by the properties of the Gamma function there exists a positive constant $K$ such that, for
any irrational number $a$
$$
\sum_{n \geq 1} n^{\gamma - 2/p }|\cos(2\pi k a)|^n  \leq   \frac{K}{(1 - | \cos (2 \pi ka) | )^{\gamma - 2/p + 1}} \, .
$$
Since $(1-|\cos(\pi u)|) \geq \pi
(d(u, {\mathbf{Z}}))^2$, we derive that
$$
\sum_{n \geq 1} n^{\gamma - 2/p }\sum_{k \in {\mathbf{Z}}^*} |k|^{-s - \varepsilon}
|\cos(2\pi k a)|^n  \leq   \frac{K}{\pi^{\gamma - 2/p + 1}} \sum_{k \in {\mathbf{Z}}^*} \frac{|k|^{-s} \big ( \log ( 1 + |k|) \big )^{-(1+\varepsilon)} }{(d(2ka, {\mathbf{Z}}) )^{2\gamma - 4/p + 2}} \, .
$$
Note that, if $a$ is badly approximable by
rationals, then so is $2a$.
Therefore if $a$ satisfies \eqref{badly}, proceeding as in the proof of Lemma 5.1 in Dedecker and Rio (2008),  we get that
\begin{equation*}
\sum_{k=2^N}^{2^{N+1}-1} \frac{1}{(d(2ka, {\mathbf{Z}}) )^{2\gamma - 4/p + 2}} \ll
2^{(2\gamma - 4/p + 2)N }  \, .
\end{equation*}
Therefore \begin{align*}
\sum_{n \geq 1} n^{\gamma - 2/p }& \sum_{k \in {\mathbf{Z}}^*} |k|^{-s} \big ( \log ( 1 + |k|) \big )^{-(1+\varepsilon)}
|\cos(2\pi k a)|^n  \\
& \ll
\sum_{N \geq 0} 2^{(2\gamma - 4/p + 2)N }  \max_{2^N \leq k \leq
2^{N+1}} |k|^{-s} \big ( \log ( 1 + |k|) \big )^{-(1+\varepsilon)}< \infty \, ,
\end{align*}
proving, by the choices of $s$ and $\gamma$, that the condition (\ref{aimcircle})  is
satisfied. This ends the proof of Corollary \ref{circle0}. Corollary \ref{circle} follows the same lines and the main point is to prove that for $\gamma$ defined by \eqref{egalitegamme} and $a$ satisfying \eqref{badlyweak},
$$
\sum_{n \geq 1} n^{\gamma - 2/p }\sum_{k \in {\mathbf{Z}}^*} |k|^{-s- \varepsilon}
|\cos(2\pi k a)|^n  < \infty \, .
$$
When $p=4$ this will give the rate $b_n=n^{1/4} (\log n )^{3/4 + \delta}$, for a $\delta >0$, in the strong invariance principle \eqref{strassen}. For $p \in ]2,4[$, the arguments also hold for some $p' >p$ leading to the rate $b_n = n^{1/p}$ in \eqref{strassen}. $\diamond $

\subsection{Proof of Corollary \ref{cortheta}.} To prove the result, we apply Corollary \ref{coralpha}, so that we need to check the conditions (\ref{Cond1cob*}), (\ref{cond2coralpha}) and (\ref{cond2coralpha*}).

Notice  that since $p >2$, (\ref{cond2coralpha}) and the first part of (\ref{Cond1cob*}) are both satisfied if
\beq \label{pourapplialpha}
\sum_{n \geq 2} \frac{n^{p-1}}{n^{2/p  } (\log n)^{\frac{(t-1)p}{2 }}}   \max_{0\leq k \leq n}\Vert \bkE_0 (X_k)\bkE_0(X_n)\Vert^{p/2}_{p/2}  < \infty \, .
\eeq We start by proving this condition.  Using Proposition 1 in Dedecker and Doukhan (2003), we get that for any integers $k$ and $n$,
$$
\Vert \bkE_0 (X_{k})\bkE_0(X_{n})\Vert_{p/2}  = \E \big ( Z_0 |\bkE_0(X_{k})| X_n\big ) \leq  \int_0^{G(\Vert \E( X_n  | {\cal F}_0 )
\Vert_1)}  Q_{|Z_0\bkE_0(X_{k})|} (u) Q_{|X_n|}(u) du \, ,
$$
where $Z_0 = \Vert \bkE_0 (X_{k})\bkE_0(X_{n})\Vert_{p/2}^{1-p/2}|\bkE_0 (X_{k})\E_0 (X_{n} )|^{p/2 -1} {\rm sign} (\E_0 (X_{n} ))$. Next using successively Lemma 2.1 and Inequality (4.6) in Rio (2000),  we obtain that
$$
\Vert \bkE_0 (X_{k})\bkE_0(X_{n})\Vert_{p/2} \leq  \int_0^{G(\Vert \E( X_n  | {\cal F}_0 )
\Vert_1)}  Q_{|Z_0|} (u) Q^2(u) du\,  .
$$
Now, applying H\"older's inequality and noting that  $\|Z_0\|_{p/(p-2)}=1$, we obtain that
$$
\max_{ 0 \leq k \leq n }\Vert \bkE_0 (X_{k})\bkE_0(X_{n})\Vert^{p/2}_{p/2} \leq   \int_0^{G (\Vert \E( X_n  | {\cal F}_0 )
\Vert_1)}  Q^p(u) du  \, ,
$$
so that (\ref{pourapplialpha}) is satisfied under (\ref{Condtheta}),
 since $ \Vert \E( X_n  | {\cal F}_0 )
\Vert_1\leq \lambda_2(n)$.

Let us prove now that the second part of (\ref{Cond1cob*}) is satisfied under (\ref{Condtheta}). Using the same arguments as before, we infer that
\begin{align*}
\Vert \bkE_0(X_{n}) \Vert_2   \leq   \Big ( \int_0^{G (\Vert \E( X_n  | {\cal F}_0 )
\Vert_1)}  Q^2(u) du \Big )^{1/2} \leq  \Big ( \int_0^1 Q^p(u) du \Big )^{1/p -1/2}  \Big (  \int_0^{G (\Vert \E( X_n  | {\cal F}_0 )
\Vert_1)}   Q^{p}(u) du \Big )^{1/2} .
\end{align*}
For the second inequality we have used H\"older's inequality together with the fact that, since $Q$ is nonincreasing, for any $x \in [0,1]$,
$$
\int_0^1 Q^p(u) du \leq \frac 1x \int_0^x Q^p (u) du \, .
$$ Therefore the second part of (\ref{Cond1cob*}) will be satisfied if
$$
\sum_{n \geq 2} \frac{n^{3p/4}}{n^{2 } (\log n)^{\frac{(t-1)p}{2 }}}
\Big (  \int_0^{\Vert \E( X_n  | {\cal F}_0 )
\Vert_1} Q^{p-1} \circ G (u) du \Big )^{p/4} < \infty\, ,
$$
which clearly holds if (\ref{Condtheta}) does, since $ \Vert \E( X_n  | {\cal F}_0 )
\Vert_1\leq \lambda_2(n)$,  $  (   \int_0^{\lambda_2(n)} Q^{p-1} \circ G (u) du  )_n$ is nonincreasing and $p>2$.

It remains to prove that (\ref{cond2coralpha*}) holds if (\ref{Condtheta}) does. Let $B_{i,j }=  X_{i}X_{j} - {\mathbf E}(X_{i}X_{j})$ and $B_{0}=|{\mathbf{E}}_{0}(B_{i,j })|^{p/2-1} \Vert {\mathbf{E}}_{0}(B_{i,j })\Vert_{p/2}^{1-p/2} {\rm sign} ( {\mathbf{E}}_{0}(B_{i,j })) $. Applying again Proposition 1 in Dedecker and Doukhan (2003), we derive that
$$
\| {\mathbf E}(X_iX_{j} | {\mathcal F}_0) - {\mathbf E}(X_{i}X_{j})\|_{p/2}
= {\mathbf{E}} ( B_{0}B_{i,j } ) \leq \int_0^{G_{B_{i,j} }(\Vert \E_0( B_{i,j} )
\Vert_1)}  Q_{B_0} (u) Q_{B_{i,j}}(u) du \, .
$$
Since $\|B_0\|_{p/(p-2)}=1$, H\"older's inequality gives
\beq \label{maj1theta}
\| {\mathbf E}(X_iX_{j} | {\mathcal F}_0) - {\mathbf E}(X_{i}X_{j})\|_{p/2} \leq  \Big ( \int_0^{G_{B_{i,j} }(\Vert \E_0( B_{i,j} )
\Vert_1)}   Q^{p/2}_{B_{i,j}}(u) du \Big )^{2/p}\, .
\eeq
Now, since $Q_{B_{i,j}}(u) \leq  Q_{|X_iX_j|}(u) + \E (|X_iX_j|)$, we get that
\begin{equation*} \label{maj1intQ}\int_0^x Q_{B_{i,j}}(u) du \leq 2\int_0^x Q_{|X_iX_j|}(u)du \leq 2 \int_0^x Q^{2}(u)du \, ,
\end{equation*}
where for the last inequality we have used Lemma 2.1 in Rio (2000). It follows that $
G_2(u/2) \leq G_{B_{i,j}} (u)
$
where  $G_2$ is the inverse of $x
\mapsto H_2(x)= \int_0^{x} Q^2 (u) du $. In particular, $G_{B_{i,j}}(u)\geq
G_2(u/c)$ for any $c \geq 2$. Since $Q$ is non-increasing, it follows that for any $c \geq 2$,
$$
\int_0^{G_{B_{i,j} }(\Vert \E_0( B_{i,j} )
\Vert_1)}   Q^{p/2}_{B_{i,j}}(u) du  =
  \int_0^{\Vert \E_0( B_{i,j} )
\Vert_1}  Q^{p/2 - 1}_{B_{i,j}}\circ G_{B_{i,j}}(u) du
  \leq \int_0^{\Vert \E_0( B_{i,j} )
\Vert_1}   Q^{p/2 - 1}_{B_{i,j}} \circ G_2(u/c) du
$$
so that
\begin{multline}\label{p1theta}
\int_0^{G_{B_{i,j} }(\Vert \E_0( B_{i,j} )
\Vert_1)}   Q^{p/2}_{B_{i,j}}(u) du \\
\leq c \int_0^{\Vert \E_0( B_{i,j} )
\Vert_1/c} Q^{p/2 - 1}_{B_{i,j}} \circ G_2(u)du
   =c \int_0^{G_2 ( \Vert \E_0( B_{i,j} )
\Vert_1/c)} Q^{p/2 - 1}_{B_{i,j}} (v) Q^2(v) dv .
\end{multline}
Notice now  that
$Q^{p/2}_{B_{i,j}}(u) \leq 2^{p/2} Q^{p/2}_{|X_iX_j|}(u) + 2^{p/2} \E (|X_iX_j|^{p/2})$, which implies that
\begin{equation} \label{majintQ}\int_0^x Q^{p/2}_{B_{i,j}}(u) du \leq 2^{1+p/2} \int_0^x Q^{p/2}_{|X_iX_j|}(u)du \leq 2^{1+p/2} \int_0^x Q^{p}(u)du \, ,
\end{equation}
where for the last inequality we have used Lemma 2.1 in Rio (2000). Therefore starting from (\ref{maj1theta}) and  using (\ref{p1theta}) together with H\"older's inequality and (\ref{majintQ}), we get that
\beq \label{maj2theta}
\| {\mathbf E}(X_iX_{j} | {\mathcal F}_0) - {\mathbf E}(X_{i}X_{j})\|^{p/2}_{p/2} \leq c \time 2^{p/2}  2^{(p^2 -4)/(2p)} \int_0^{G_2 ( \Vert \E_0( B_{i,j} )
\Vert_1/c)}  Q^p(v) dv \, .
\eeq
We show now that for  any $i \geq n$ and $j \geq n$, there exists $c \geq 2$ such that
\beq \label{objtheta}
G_2 ( \Vert \E_0( B_{i,j} )
\Vert_1/c) \leq G (\theta_2 (n) ) \, .
\eeq
With this aim, let $M = Q\circ G (\theta_2 (n)) $ and, for $i  \in {\mathbf Z} $, define the variables
$$
X_i'= X_i \I_{| X_i
|\leq M }
 \quad \text{and} \quad    X_i''
= X_i  \I_{| X_i | > M } \, .
$$
With this notation and using the stationarity, we have that
\begin{align}  \label{maj3theta}
\Vert \E_0( B_{i,j} )
\Vert_1 & \leq \Vert \E_0(X_i'X_j' ) - \E(X_i'X_j')
\Vert_1 +4 M \E |X_0''| + 2 \E |X_0''|^2 \nonumber \\
& \leq \Vert \E_0(X_i'X_j' ) - \E(X_i'X_j')
\Vert_1 +6 \E ( X_0^2 \I_{| X_0|> M } ) \, .
\end{align}
Notice now that for any reals $x$ and $y$:
\begin{eqnarray} \label{inega1}
& & \Big | \vert x\I_{| x | \leq M } + y \I_{| y |
\leq M } \vert^2- \vert x+y  \vert^2 \I_{| x+ y |
\leq 2 M }
\Big | \nonumber \\
& &  \leq  | y |^2 \I_{| y | \leq M } \I_{| x | > M
} + | x |^2 \I_{| x | \leq M } \I_{| y |
> M } + | x+ y |^2 \I_{| x+ y | \leq 2 M } \I_{| x | > M
} \nonumber \\
& & \quad + | x+ y |^2 \I_{| x+ y | \leq 2 M } \I_{| y
| > M }+ | x+ y |^2 \I_{| x+ y | \leq 2 M } \I_{| x
| > M } \I_{| y | > M } \nonumber \\
& & \leq 5M^2 \I_{| x | > M } + 9 M^2 \I_{| y | > M
}\, .
\end{eqnarray}
In addition setting for any real $ u \geq 0$ and any $T>0$, $g_{T}
(u) = u^2\wedge T^2 $, we have
\beq \label{inegal2}
 \Big | \vert x+y  \vert^2 \I_{| x+ y |
\leq 2 M } - g_{2M} \big( \vert x+y  \vert \big ) \Big | \leq
4M^2 \I_{\vert x+y  \vert > 2 M } \, .
\eeq
Using (\ref{inega1}) and (\ref{inegal2}), the fact that $
4xy= | x + y|^2 - | x - y|^2 $, and the stationarity  of $(X_i)_{i \in {\mathbf Z}}$, it follows that
\begin{equation} \label{inegal3}
\sup_{j\geq k \geq q }  \Vert 4 X'_jX_k'-  g_{2M} \big(
\vert X_j+X_k   \vert \big )   +  g_{2M} \big( \vert X_j - X_k
\vert\big ) \Vert_1
 \leq  44  M \E (| X_0| \I_{| X_0|> M } )\leq 44  \E ( X_0^2 \I_{| X_0|> M } )\,
.
\end{equation}
In addition, since $g_{2M}$ is $2M$-Lipschitz, it
follows that
\beq \label{inegal4}
 \sup_{j\geq k \geq n } \Vert \E \big ( g_{2M} (
\vert X_j+X_k  \vert) | {\cal F}_{0} \big ) -   \E \big (
g_{2M} ( \vert X_j+X_k \vert)  \big ) \Vert_1 \leq 2M \theta_2
(n) \, ,
\eeq
and the same holds true with $\vert X_j+X_k  \vert$ in place
of $\vert X_j-X_k  \vert $. Starting from (\ref{maj3theta}) and using (\ref{inegal3}), (\ref{inegal4}), together with the fact that $\E ( X_0^2 \I_{| X_0|> M } ) \leq \int_0^{G ( \theta_2 (n))} Q^2 (u) du $, we obtain that for any $i \geq n$ and $j \geq n$,
$$
\Vert \E_0( B_{i,j} )
\Vert_1  \leq Q\circ G (\theta_2 (n))   \theta_2
(n) +  28 \int_0^{G ( \theta_2 (n))} Q^2 (u) du  \leq   29 \int_0^{ \theta_2 (n)} Q \circ G (u) du
$$
so that
$$
\Vert \E_0( B_{i,j} )
\Vert_1 \leq 29 \int_0^{G ( \theta_2 (n))}  Q^2 (u) du = 29  H_2 ( G (\theta_2 (n)) ) \, ,
$$
which proves (\ref{objtheta}) with $c =29$.
Starting from (\ref{maj2theta}),  using (\ref{objtheta}) and the fact that $\theta_2(n) \leq \lambda_2(n)$, it follows that  (\ref{cond2coralpha*}) will be satisfied if
$$
\sum_{n \geq 2} \frac{n^{p-2}}{ (\log n)^{\frac{(t-1)p}{2 }}}  \int_0^{\lambda_2(n)} Q^{p-1} \circ G (u) du  < \infty\, ,
$$
which holds if (\ref{Condtheta}) does, since $  (   \int_0^{\lambda_2(n)} Q^{p-1} \circ G (u) du  )_n$ is nonincreasing and $p>2$. This ends the proof of Corollary \ref{cortheta}. $\diamond$

\subsection{Proof of Corollary \ref{coralphafaible}.} It comes from an application of Corollary \ref{coralpha}. We omit the details since to prove that  (\ref{Cond1cob*}), (\ref{cond2coralpha}) and (\ref{cond2coralpha*}) are  satisfied if (\ref{Condstrong}) is, we follow the lines of the beginning of the proof of Corollary \ref{cortheta} and, to take care of the covariances terms, we use the arguments developped in the proof of Proposition 5.3 in Merlev\`ede and Rio (2012). $\diamond$

\section{Appendix}
\label{sectionappendix}
\setcounter{equation}{0}
The next proposition gives simple criteria to obtain rates of convergence in the almost sure invariance principle for a strictly stationary sequence of martingale differences. It is based on Theorem 2.1 in Shao (1993) that gives, with the help of the Skorohod-Strassen theorem, sufficient conditions for partial sums of non necessarily stationary sequences of martingale differences to satisfy the almost sure invariance principle with rates.
\begin{Proposition} \label{propmart}  Let $d_i=d_0 \circ T^i$ where $d_0$ is an element of $H_0 \ominus H_{-1}$. Assume furthermore that $d_0$ is in ${\mathbf L}^p$ with $p\in ]2,4]$. Let $M_n = \sum_{i=1}^n d_i$.  Let $\psi$ be a nondecreasing positive function. Assume also that there exists a positive constant $C$ such that $\psi(2n) \leq C \psi(n) $ for all $n \geq 1$ and $\sum_{n>0}
 \psi^{-p/2}(n)
<\infty$. Define $v_n=n^{-1}\psi^{p/2}(n)\sum_{k\ge n}\psi^{-p/2}(k)$. Assume that\beq \label{condcarremart}
 \sum_{n \geq 1} v_n^{p/2} \frac{\big \Vert  \bkE (M_n^2 | {\mathcal F}_0  ) -\bkE  (M_n^2)
 \big \Vert^{p/2}_{p/2} }{ n \psi^{p/2} (n) }< \infty \, .\eeq Enlarging $\Omega$ if
necessary, there exists a sequence $(Z_i)_{i
\geq 1}$ of iid gaussian random variables with zero mean and
variance ${\mathbf E}(d_0^2)$ such that
$$
\sup_{1\leq k \leq n} \Big|M_k - \sum_{i=1}^k Z_i\Big|  = o \Big ( \psi^{1/2}(n) \Big |  \log \frac{ n}{\psi(n)} + \log \log (\psi(n)) \Big |^{1/2}  \Big ) \text{ almost surely} \, .
$$
\end{Proposition}

\noindent {\bf Proof of Proposition \ref{propmart}.}
From Theorem 2.1 in Shao (1993), it suffices to verify that
\beq \label{but1}
 \sum_{i=1}^n (  {\mathbf E}(d^2_i|{\mathcal F}_{i-1}) -{\mathbf E}(d_i^2))  = o(\psi(n)) \text{ almost surely}\, .
\eeq
Since $\psi$ is nondecreasing and 
$\psi(2n)\le C\psi(n)$ for any
$n\ge 1$  and  some fixed $C>0$, the almost sure convergence \eqref{but1} will follow if we can prove that,
for every $\varepsilon>0$, 
$$\sum_{r\geq1}\mathbf{P}\Big(\max_{1\leq k\leq 2^{r+1}}\Big|\sum_{i=1}^k {\mathbf E}(d^2_i|{\mathcal F}_{j-1}) -{\mathbf E}(d_i^2))\Big|\geq
\varepsilon \psi(2^{r})\Big)<\infty\, .$$ Therefore \eqref{but1} and then the proposition will be proved if we can show that
$$
\sum_{r\geq1}\frac {1}{\psi^{p/2}(2^r)}{\E \Big( \max_{1\leq k\leq 2^{r}}\Big |\sum_{i=1}^k {\mathbf E}(d^2_i|{\mathcal F}_{i-1}) -{\mathbf E}(d_i^2))\Big |^{p/2}\Big )} <\infty\, .
$$
Applying Theorem 3 in Wu and Zhao (2008) (since $1 < p/2 \leq 2$) and using the martingale property, we get that
$$
\bkE \Big (  \max_{ 1 \leq k \leq 2^r} \Big |  \sum_{i=1}^k (  {\mathbf E}(d^2_i|{\mathcal F}_{i-1}) - {\mathbf E}(d^2_i))  \Big |^{p/2}  \Big ) \ll  2^r \| d^2_1 \|^{p/2}_{p/2}  + 2^r \left ( \sum_{k=0}^{r-1} \frac{\|{\mathbf E}(M_{2^k}^2|{\mathcal F}_{0}) - {\mathbf E}(M_{2^k}^2) \|_{p/2}}{2^{2k/p}}  \right )^{{p/2}} \, .
$$
Since $\psi$ is nondecreasing and since $\sum_{n} \psi^{-p}(n)<\infty$, it follows that
$\sum_r 2^r\psi^{-p}(2^r)<\infty$. Hence the result will be proved if we can show that
$$
\sum_{r\ge 1} \frac{1}{\psi^{p/2}(2^r)} \Big ( \sum_{k=0}^{r-1} 2^{2(r-k)/p} \|{\mathbf E}(M_{2^k}^2|{\mathcal F}_{0}) - {\mathbf E}(M_{2^k}^2) \|_{p/2} \Big )^{{p/2}}<\infty \, .
$$
But the latter is implied by \eqref{condcarremart}, which may be proved as in the proof of $(3)\Rightarrow (2)$ of Proposition 2.2 in Cuny (2011) (see Appendix A there).
$\diamond$

\medskip

The following lemma is useful to compare conditions involving the projection operator $P_0$ with a mixingale-type condition.
\begin{Lemma} \label{complemmap0xi} Let $p\geq 2$. For any real $1 \leq q \leq p$ and any positive integer $n$,
$$
\sum_{k \geq 2 n}\Vert P_0(X_{k}) \Vert^q_p \ll \sum_{k \geq n  } \frac{\Vert \E(X_k |{\mathcal F}_0 )\Vert^q_p}{k^{q/p}} \, .
$$
\end{Lemma}

\noindent {\bf Proof of Lemma \ref{complemmap0xi}.}
If the upper bound is infinite, the inequality is clear. 

Let us consider now the case where the upper bound is finite. In that case, since 
$\|\E(X_n |{\mathcal F}_0 )\|_p$ is nonincreasing, we infer that
$n^{1-q/p}\|\E(X_n |{\mathcal F}_0 )\|_p^q$ converges to $0$ as $n$ tends
to infinity, so that ${\mathbb E}(X_0|{\mathcal F}_{-\infty})=0$ almost surely. 

Now, for any sequence of nonnegative numbers $(u_k)_{k \in {\mathbf N}}$ and any real $\alpha >1$, the following inequality holds: for any positive integer $n$, there exists a positive constant $C_\alpha$ depending only on $\alpha$ such that
\beq \label{factseries}
\sum_{k \geq 2n} u_k \leq C_\alpha \sum_{k\geq n+1} \frac{1}{k^{1/\alpha}} \Big ( \sum_{\ell  \geq k } u_{\ell}^\alpha \Big )^{1/\alpha} 
\eeq
(to prove this inequality, it suffices to slightly adapt the proof of Lemma 6.1 in Dedecker {\it et al.} (2011) and to use the fact that the sequence $\big ( \sum_{\ell  \geq k } u_{\ell}^\alpha \big )_{k \geq 1}$ is nonincreasing).  

We proceed now as in the proof of Remark 3.3 in  Dedecker {\it et al.} (2011). We first consider the case  $p>q$. Applying inequality \eqref{factseries} with $u_k = \Vert P_0(X_{k}) \Vert^{q}_p$ and $\alpha=p/q$, we get that
\begin{equation} \label{ros0}
\sum_{k \geq 2 n}\Vert P_0(X_{k}) \Vert_p^q  \ll \sum_{k\geq n+1} \frac{1}{k^{q/p}} \Big ( \sum_{\ell  \geq k } \Vert P_0(X_{\ell}) \Vert^p_p \Big )^{q/p}
\end{equation}
Next, using the stationarity and applying the Rosenthal's inequality given in Theorem 2.12 of Hall and Heyde
(1980), we then derive that for any $p\in \lbrack 2,\infty \lbrack $, there
exists a constant $c_{p}$ depending only on $p$ such that
\begin{equation} \label{ros1}
 \sum_{\ell  \geq k } \Vert P_0(X_{\ell}) \Vert^p_p =  \sum_{\ell  \geq k } \Vert P_{-\ell}(X_{0}) \Vert^p_p\leq c_{p}\Big \Vert
\sum_{\ell  \geq k } P_{-\ell}(X_{0}) \Big \Vert _{p}^{p}=c_{p}\Vert {\mathbf{E}}%
(X _{k}|{\mathcal{F}}_{0})\Vert _{p}^{p}\,,
\end{equation}%
 the last equality being true because $\mathbf E(X_0 |{\mathcal F}_{-\infty}) =0$ 
 almost surely. Therefore when $p>q$, the lemma follows by taking into account \eqref{ros1} in \eqref{ros0}. Now when $p=q$, inequality \eqref{ros1} together with the fact that $\Vert {\mathbf{E}}(
X _{k}|{\mathcal{F}}_{0})\Vert _{p}$ is nonincreasing implies the result. $\diamond$

\bigskip

The next lemma is useful to deal with sequences having a subadditive property (see Lemma 38 in Merlev\`ede and Peligrad (2012) for a proof).

\begin{Lemma}
\label{claim1} Let $(V_{i})_{i \geq0}$ be a sequence of non negative
numbers such that $V_{0} = 0$ and for all $i,j \geq0$,
\begin{equation*}
\label{condsubadd}V_{i+j} \leq C ( V_{i} + V_{j}) \, ,
\end{equation*}
where $C \geq1$ is a constant not depending on $i$ and $j$. Then
 For any integer $r \geq1$, any integer $n$ satisfying $2^{r-1} \leq n <
2^{r}$ and any real $q \geq0$
\[
\sum_{i=0}^{r-1} \frac{1}{2^{iq}} V_{2^{i}} \leq C 2^{q+2} ( 2^{q+1} -1)^{-1}
\sum_{k= 1}^{n} \frac{V_{k}}{k^{1+q}} \, .
\]
\end{Lemma}

\bigskip

\noindent\textbf{Acknowledgment}. The authors are indebted to the referees for
carefully reading the manuscript and for their  helpful comments. They also wish to thank Sergey Utev for helpful discussions.

\end{document}